\documentclass{article}

\usepackage[english]{babel}
\usepackage{amssymb}
\usepackage{mathtools}
\usepackage{graphicx}
\usepackage{caption}
\usepackage{subcaption}
\usepackage{epsfig}
\usepackage{amsmath}
\usepackage{amsthm}
\usepackage{setspace}
\usepackage{arxiv}
\usepackage{mathrsfs}
\usepackage{enumitem}
\usepackage{subcaption}
\usepackage{color}
\usepackage{tcolorbox}
\usepackage[ruled,vlined]{algorithm2e}
\SetKwProg{initialize}{Initialize:}{}{}
\SetKw{KwBy}{by}
\usepackage{gensymb}
\usepackage{multirow}
\usepackage{adjustbox}
\usepackage[graphicx]{realboxes}
\usepackage{longtable}
\usepackage{hhline}
\usepackage{geometry}
\usepackage[
  backend=bibtex,
  natbib=true,
  sorting=none,
  maxcitenames=2,
  mincitenames=1,
  useprefix=true
]{biblatex}

\usepackage{booktabs}
\usepackage{lipsum}                     
\usepackage{xargs}                      
\usepackage[ruled,vlined]{algorithm2e}

\usepackage[colorinlistoftodos,prependcaption,textsize=tiny]{todonotes}
\newcommandx{\unsure}[2][1=]{\todo[linecolor=red,backgroundcolor=red!25,bordercolor=red,#1]{#2}}
\newcommandx{\change}[2][1=]{\todo[linecolor=blue,backgroundcolor=blue!25,bordercolor=blue,#1]{#2}}
\newcommandx{\info}[2][1=]{\todo[linecolor=OliveGreen,backgroundcolor=OliveGreen!25,bordercolor=OliveGreen,#1]{#2}}
\newcommandx{\improvement}[2][1=]{\todo[linecolor=Plum,backgroundcolor=Plum!25,bordercolor=Plum,#1]{#2}}
\newcommandx{\thiswillnotshow}[2][1=]{\todo[disable,#1]{#2}}
\theoremstyle{plain}
\newtheorem{thm}{Theorem}[section]
\newtheorem{lem}[thm]{Lemma}
\newtheorem{prop}[thm]{Proposition}
\newtheorem{cor}[thm]{Corollary}

\theoremstyle{definition}
\newtheorem{definition}[thm]{Definition}

\theoremstyle{remark}
\newtheorem{remark}[thm]{Remark}
\definecolor{ForestGreen}{RGB}{34,139,34}
 
\usepackage{float}
\restylefloat{table}
\DeclareMathOperator*{\argmin}{arg\,min}
\usepackage{bbm}

\newcounter{example}[section]

\usepackage{nomencl}
\makenomenclature

\usepackage{multicol}
\usepackage{color}
\usepackage[hidelinks]{hyperref}

\usepackage{comment}
\title{Controlling avalanche sizes in the Abelian sandpile model}
\date{\today}

\author{
\begin{tabular}{cc}
Maike C. de Jongh\textsuperscript{a}
&
Antal A. Járai\textsuperscript{b}
\\
\texttt{m.c.dejongh@utwente.nl}
&
\texttt{a.jarai@bath.ac.uk}
\\[0.35cm]
Richard J. Boucherie\textsuperscript{a}
&
M.N.M. van Lieshout\textsuperscript{c,a}
\\
\texttt{r.j.boucherie@utwente.nl}
&
\texttt{m.n.m.van.lieshout@cwi.nl}
\end{tabular}
\\[1.5cm]
\textsuperscript{a}Department of Applied Mathematics, University of Twente,\\
P.O. Box 217, 7500 AE Enschede, The Netherlands
\\[0.15cm]
\textsuperscript{b}Department of Mathematical Sciences, University of Bath,\\
Claverton Down, Bath BA2 7AY, United Kingdom
\\[0.15cm]
\textsuperscript{c}Centrum Wiskunde \& Informatica,\\
P.O. Box 94079, 1090 GB Amsterdam, The Netherlands
}

\begin{document}
\raggedbottom
\newcommand{\edt}[1]{{\vbox{ \hbox{#1} \vskip-0.3em \hrule}}}
\newcommand{\RR}{{\mathbb R}}
\newcommand{\QQ}{{\mathbb Q}}
\newcommand{\NN}{{\mathbb N}}
\newcommand{\cf}[1]{{\bf{#1}}}
\newcommand{\EE}{{\mathbb E}}
\newcommand{\PP}{{\mathbb P}}
\newcommand{\cA}{{\mathcal A}}
\newcommand{\cX}{{\mathcal X}}
\newcommand{\cH}{{\mathcal H}}

\newcommand{\expect}{\operatorname{\EE}\expectarg}
\DeclarePairedDelimiterX{\expectarg}[1]{[}{]}{%
  \ifnum\currentgrouptype=16 \else\begingroup\fi
  \activatebar#1
  \ifnum\currentgrouptype=16 \else\endgroup\fi
}
\DeclarePairedDelimiter\ceil{\lceil}{\rceil}
\DeclarePairedDelimiter\floor{\lfloor}{\rfloor}

\newcommand{\innermid}{\nonscript\;\delimsize\vert\nonscript\;}
\newcommand{\activatebar}{%
  \begingroup\lccode`\~=`\|
  \lowercase{\endgroup\let~}\innermid 
  \mathcode`|=\string"8000
}

\hyphenation{high-di-men-sio-nal}

\allowdisplaybreaks

\maketitle

\begin{abstract}
    The Abelian sandpile model is a classical framework for studying large cascading events, such as those arising in models of forest fires and earthquakes. This work studies the control of such avalanches through the targeted removal of sand grains. We formulate this sequential control problem as an average-cost Markov decision process in which, before each grain addition, a controller may remove a sand grain from a selected vertex. The objective is to balance the long-run expected avalanche size against the cost of intervention.  We analyze local intervention policies, which remove grains from a prescribed set of vertices whenever possible, and establish structural results for the recurrent classes and stationary distributions of the induced Markov chains. When the intervention set satisfies a property called topple-freeness, we sharpen these results, derive expressions for the gain and bias functions, and give necessary and sufficient conditions for optimality on general finite graphs. Finally, we apply the general theory to several specific graph classes and obtain a complete characterization of optimal policies on the one-dimensional lattice, the complete graph and the star graph for every value of the intervention cost.  
\end{abstract}

\setcounter{tocdepth}{3}
\tableofcontents

\section{Introduction}

The Abelian sandpile model, introduced by \citet{Bak1, Bak2, Dhar}, provides a canonical framework for studying self-organized criticality. The model describes a dynamical system defined on a graph in which sand grains are successively added at randomly selected vertices. When the number of grains at a vertex exceeds its capacity, the vertex topples and distributes grains among its neighbors. This may cause topplings of other vertices, giving rise to an avalanche that can propagate across a large part of the graph. Grains 
leave the system at a sink vertex, playing the role of boundary. Critical behavior emerges in the model despite the absence of a parameter that can be tuned to drive the system towards a critical point. This mechanism, known as self-organized criticality (SOC) \cite{Jensen_1998}, has been associated with phenomena such as earthquakes \citep{Scholz}, forest fires \citep{Drossel}, fluctuations in financial markets \citep{Bartolozzi}, and evolutionary dynamics \citep{Bak3}. 

Since large-scale cascades 
can have destructive consequences in applied settings, a natural question is to what extent 
their occurrence or impact can be reduced by means of external control. 
In this paper, we use the framework of Markov decision processes (MDP) \cite{Puterman} to give an analysis of the above question in the context of the Abelian sandpile model, focusing on strategies that are optimal in the long-run. This extends the earlier work of \citet{deJongh1}, which studied the immediate rather than the long-run effect of interventions. While the Abelian sandpile model only shares broad qualitative features with real avalanches, we believe that its study provides useful insight into controlled SOC. The main difficulty in the analysis of SOC is the inherent strong spatial and temporal correlations of the cascades, and some aspects of these features are reasonably well understood at the mathematical level in the sandpile model. Namely, when the cost associated with an avalanche is \emph{linear} in the number of topplings involved, Dhar's formula \cite{Dhar} relates expected toppling numbers to the Green's function of the underlying graph. This relation plays an important role in our proof of Theorem \ref{gain_bias_never_intervene} and in the explicit calculations for specific graph classes in Section \ref{Applications_to_specific_graph_classes}.

The setup of our work is as follows. We formulate the problem of controlling avalanche sizes through grain-removal interventions as an average-cost MDP. Before each grain addition, the controller may either take no action or remove one grain from a selected vertex. The cost function of the MDP is the sum of an immediate intervention cost $\beta$ (when the controller decides to intervene) and the expected size of the subsequent avalanche, thereby capturing a trade-off. Since new sand grains arrive at each time step, our aim is to identify intervention strategies that achieve a sustained reduction in expected avalanche size relative to their intervention costs. The objective, therefore, is to minimize the long-run average cost per grain addition, or the \textit{gain}. This formulation also allows us to exploit the stationary structure of the controlled sandpile dynamics to derive explicit expressions for the gain and bias and establish conditions for optimality. Greater weight could instead be assigned to costs in the near future by using a discounted cost criterion~\cite{Puterman}. This provides a different perspective that we do not consider in this work.


We study \emph{local intervention policies} (see Definition \ref{local_intervention}), which remove grains from a prescribed set $W$
of vertices whenever possible. We establish structural results for the recurrent classes and stationary distributions of the Markov chains induced by these policies. The special case $W = \emptyset$ corresponds to the \emph{never-intervene policy}. 
We obtain stronger results when the intervention set $W$  is \emph{topple-free} under the policy, which means that in stationarity, no vertex in $W$ topples during a transition. We give the precise definition in Definition \ref{topple_free} 
and identify several graph classes and intervention sets for which this property holds.  Under the topple-free condition, we obtain sharper characterizations of recurrent classes and stationary distributions, as well as explicit expressions for the gain (i.e., the long-run average cost) and bias (i.e., the relative transient cost) functions (precise definitions of these are given in Section \ref{ssec:optimality-eqns}) and conditions for the optimality of the corresponding policy. 
In addition to these results being of relevance for the control problem, we believe they are also of interest from the point of view of the theory of Abelian sandpiles:  under the never-intervene policy, we recover the usual sandpile dynamics, whereas the controlled Markov chains represent perturbations of the usual sandpile dynamics that preserve some, but need not retain all, of the structure of the original model.

One of our main contributions is the following theorem, stated here somewhat informally;
see Theorem \ref{general_optimal_policies} for the precise version.
\begin{thm}\label{main_gen}
Consider a local intervention policy with intervention set $W$, and suppose that 
$W$ is topple-free under this policy. There are explicit bounds 
$\underline{\beta}_W, \overline{\beta}_W \in \RR$, depending only on $W$ and the graph, such that if $\underline{\beta}_W \le \overline{\beta}_W$, then for every $\beta \in [\underline{\beta}_W, \overline{\beta}_W]$
the policy is optimal.
\end{thm}
\noindent
The thresholds $\underline{\beta}_W, \overline{\beta}_W$ are expressed in terms of random 
walk quantities on the underlying graph.

As an application, we prove the following theorem, that identifies optimal policies 
for every $\beta \ge 0$ in one dimension. Consider the one-dimensional lattice $\{1, \ldots, L\}$ with vertices $1$ and $L$ connected to the sink by a single edge. Let $K = \lfloor (L-1)/2 \rfloor$ and $W_k = \{k+1, \ldots, L-k\}$ for $k = 0, \ldots, K$. Define
\[
\beta_k = \dfrac{k(k+1)}{2}, \qquad k = 0, \ldots, K, \qquad \beta_{K+1} = \begin{cases}
    (L+1)^2/8, &\text{if } L \text{ is odd,}\\
    L(L+2)/8, &\text{if } L \text{is even.}
\end{cases}
\]
Let $\pi_0$ be the always-intervene policy that removes the leftmost grain whenever the configuration is nonempty and for each $k = 1, \ldots, K$, let $\pi_k$ be any local intervention policy with intervention set $W_k$.
\begin{thm}\label{main_1D}
For each $k = 0, \ldots, K$, the policy $\pi_k$ is optimal if and only if $\beta \in [\beta_k, \beta_{k+1}]$. The never-intervene policy is optimal if and only if $\beta \geq \beta_{K+1}$. 
\end{thm}

Our further contributions are the following. We derive general expressions for the gain and bias of the never-intervene policy and an exact threshold for its optimality on finite graphs. These results provide the foundation for the analysis of local intervention policies with general intervention sets. We apply the general theory for several graph classes, obtaining explicit quantities and graph-specific optimality results. For the complete graph and the star graph, these results yield a full characterization: for every value of $\beta$, an optimal policy can be found among the policies that either never intervene or always intervene. 

\subsection{Related work}
\label{ssec:related-work}
Several intervention mechanisms have been studied for SOC systems. One line of research seeks to influence avalanche sizes by selecting the addition location of the next grain based on the current state of the system \citep{Falk, Noel, Parsaeifard}. Related strategies deliberately trigger avalanches at vertices that are close to becoming unstable, aiming to prevent large-scale avalanches through more frequent, localized dissipation events \citep{Cajueiro1, Cajueiro2}. Since the addition site cannot be controlled in many applications, \citet{Qi} instead introduce failure tolerance, which temporarily prevents the toppling of overloaded vertices. In a different model of cascading failure, \citet{Hou} investigate how heterogeneous redistribution of the load of such vertices affects the size of the resulting avalanche. Control mechanisms for cascading failures have also been studied in interdependent networks \citep{Turalska}. In a model of systemic risk in financial networks, \citet{Scala} study targeted immunization policies that provide assistance to selected institutions. 

The studies above rely primarily on numerical methods to investigate strategies for reducing avalanche sizes. A rigorous analysis of interventions was provided by \Citet{deJongh1}, who consider removing sand grains from square-shaped connected regions consisting of vertices at maximal height. They quantify how such interventions affect the size of the next avalanche and identify optimal intervention locations that minimize its expected size. This paper extends their analysis from the effect of a single intervention on the next avalanche to the long-term performance of intervention strategies. We formulate this sequential control problem as an average-cost Markov decision process (MDP)~\citep{Puterman}.

MDPs have proved useful in a variety of decision-making problems involving spatially interacting or networked systems. For example, the framework has been used to derive optimal strategies for controlling droplet growth in the Ising model \citep{deJongh2, deJongh3, Baldassarri}. Other applications include the containment and suppression of wildfires \citep{Haksar2018, Haksar2019}, epidemic control on social networks and under uncertain transmission dynamics \citep{Nasir2023, Song2023}, the control of colloidal self-assembly processes for producing crystalline structures \citep{Tang2016, Tang2017} and optimal thinning strategies in forest management \citep{Lieshout}. 

A related sequential control problem was studied by \citet{Cajueiro3}, who introduce a dynamic programming formulation for the control of the directed Abelian Dhar-Ramaswamy model \citep{DharRamaswamy1989}. Before each addition of a new sand grain, the controller may either take no action or select a vertex at which to trigger an avalanche. The objective is to balance the cost of intervention against a penalty associated with avalanche size. Due to the exponential growth of the state space, an optimal strategy is computed numerically only for a $4 \times 4$ lattice. The resulting solution is then used as a benchmark for heuristics that can be applied to larger systems. 
Although our formulation shares a sequential decision-making perspective, the two control frameworks differ in several fundamental respects. While \citet{Cajueiro3} consider discounted control of a directed sandpile model through deliberately triggered avalanches, we consider average-cost control of the standard Abelian sandpile model through grain removal. Moreover, whereas their analysis of optimal policies is numerical and restricted to a fixed lattice size, our analysis yields rigorous structural results for general finite graphs, together with more detailed results for several specific graph classes.

\subsection{Organization}
The remainder of this paper is organized as follows. Section \ref{Model definition} formulates the control problem and provides the necessary background on the Abelian sandpile model and average-cost MDPs.  Section \ref{General_results_for_local_intervention_policies} presents general results for local intervention policies, including sharper results under the topple-free condition. In Section \ref{Applications_to_specific_graph_classes}, we apply the general theory to several specific graph classes, providing a complete characterization of the optimal policies on the one-dimensional lattice, the complete graph and the star graph for every value of $\beta$. In Section \ref{Discussion}, we give some concluding remarks and directions for future research.

\section{Model definition}
\label{Model definition}
In this section, we formally define the control problem studied in this work and provide the necessary background on the Abelian sandpile model and average-cost MDPs.

\subsection{The Abelian sandpile model}\label{ASM_model}
We consider the Abelian sandpile model on a finite, undirected graph $\mathcal{G} = (V \cup \{s\}, E)$, with vertex set $V$, edge set $E$, and a distinguished vertex $s$ called the sink. Let $\deg(x)$, $x \in V$, denote the degree of $x$ in $\mathcal{G}$, including edges incident to the sink. We allow parallel edges, but assume that the graph contains no self-loops. Let $b_{xy}$, $x, y \in V$, denote the number of edges between $x$ and $y$. A sandpile configuration on $\mathcal{G}$ is a function $\eta: V \rightarrow \mathbb{Z}_{\geq 0}$, where $\eta(x)$ denotes the height of the vertex $x \in V$, that is, the number of indistinguishable sand grains at $x$. We denote the set of all sandpile configurations by $\mathcal{X}$. The total mass $m(\eta)$ of a sandpile $\eta \in \mathcal{X}$ is defined by
\[
m(\eta) = \sum\limits_{x \in V} \eta(x), \qquad \eta \in \mathcal{X}. 
\]
Let
\[
V^+(\eta) = \{x \in V: \eta(x) > 0\}
\]
denote the set of vertices containing at least one sand grain in the configuration $\eta \in \mathcal{X}$.
A vertex $x \in V$ is stable in $\eta \in \mathcal{X}$ if $\eta(x) < \deg(x)$. A configuration is stable if each vertex is stable. Let the set of all stable sandpile configurations be denoted by $\Omega$. We call a vertex $x \in V$ critical in $\eta \in \Omega$ if $\eta(x) = \deg(x)-1$. A vertex $x \in V$ is noncritical if $\eta(x) < \deg(x)-1$. For $i = 0, \ldots, \sum\limits_{x \in V} (\deg(x) - 1)$, let
\begin{equation}\label{Omega_i}
\Omega_i = \{\eta \in \Omega: m(\eta) = i\}
\end{equation}
denote the set of stable configurations containing exactly $i$ sand grains. 

The sandpile Markov chain $(\eta_t)_{t \geq 0}$ on the state space $\Omega$ is defined by the following dynamics. At each discrete time $t = 0, 1, \ldots$, a vertex $x \in V$ is selected uniformly at random and a new sand grain is added at $x$. We denote the addition operator by $\alpha_x: \mathcal{X} \rightarrow \mathcal{X}$, where $\alpha_x \eta$ is obtained from $\eta$ by adding one sand grain at $x$. If $x$ becomes unstable due to this addition, then it topples, meaning that it sends one sand grain to each of its neighboring vertices. Any vertices that become unstable subsequently also topple, and this process continues until the configuration is stable. We refer to this procedure as the relaxation of the configuration. 
Sand grains that reach the sink during the relaxation are dissipated.
Given a configuration $\eta \in \mathcal{X}$, we denote its stabilization by $\eta^{\circ}$. A fundamental result due to Dhar~\citep{Dhar}, called the Abelian property, states that, starting from a configuration $\eta \in \Omega$ and adding a grain at a vertex $x$, the number of times each vertex $y \in V$ topples during the subsequent relaxation, denoted by $n_y(\eta, x)$, is independent of the order in which unstable vertices are toppled. In particular, $(\alpha_x \eta)^{\circ}$ is well-defined.

Let $\Delta$ denote the restriction of the graph Laplacian of $\mathcal{G}$ to $V \times V$, i.e.,
\[
\Delta_{xy} = \begin{cases}
    \deg(x), &\text{if } x = y \in V; \\
    -b_{xy}, &\text{if } x \neq y, \quad x,y \in V.
\end{cases}
\]
Identifying a configuration $\eta$ with a row vector, the procedure of toppling a vertex $x$ can be written as
\[
\eta \rightarrow \eta - \Delta_{x, \cdot},
\]
where $\Delta_{x, \cdot}$ denotes the row of $\Delta$ corresponding to $x$.
The transition kernel of the sandpile Markov chain is given by
\begin{equation}\label{trans_probs_uncontrolled}
p(\eta, \eta') = \dfrac{1}{|V|}\sum\limits_{x \in V}\mathbbm{1}[(\alpha_x \eta)^{\circ} = \eta'], \qquad \eta, \eta' \in \Omega.
\end{equation}
The sandpile Markov chain has a unique recurrent class and its stationary distribution, denoted by $\mu$, is the uniform distribution on this class \cite{Dhar}. 

The main quantity of interest in this paper is the avalanche size, defined as the total number of topplings during the relaxation process. Our objective is to control the sandpile dynamics in such a way as to minimize the expected avalanche size. The avalanche size resulting from adding a grain at $x \in V$ to a configuration $\eta \in \Omega$ is therefore given by
\[
n(\eta, x) = \sum\limits_{y \in V} n_y(\eta, x).
\]
Let $N(\eta)$ denote the size of the avalanche triggered by the next grain addition starting from a configuration $\eta \in \Omega$. Note that this is a random variable with distribution
\[
\mathbb{P}(N(\eta) = k) = \dfrac{1}{|V|}\sum\limits_{x \in V} \mathbbm{1}[n(\eta, x) = k], \qquad k = 0, 1, \ldots.
\]
When we need to consider the full time-evolution $(\eta_t)_{t \ge 1}$, we distinguish the avalanche sizes as $N_1(\eta_1), N_2(\eta_2), \dots$. Here, conditional on $\eta_1 = \eta$, 
the random variable $N_1$ has the distribution of $N(\eta)$. The next configuration
$\eta_2$ is a function of $\eta_1$ and the toppling numbers 
$( n_y(\eta_1,X_1) )_{y \in V}$, where $X_1$ is the random location of the added grain.
In general, denote the history of the process up to time $t$ by
$\mathcal{F}_t = \sigma(\eta_1,X_1,\dots,\eta_t,X_t)$. Then $\eta_{t+1}$ is determined by $\eta_t, X_t$. On the event $\{ \eta_{t+1} = \eta \}$, the random variable 
$N_{t+1}(\eta_{t+1})$ is a copy of $N(\eta)$ independent of $\mathcal{F}_t$.

\subsection{The sandpile Markov decision process}
\label{ssec:sandpile-MDP}
We consider a setting in which the dynamics of the sandpile Markov chain can be controlled by means of local interventions. Specifically, before each new grain addition, an external decision maker may choose a vertex $x \in V$ from which to remove a sand grain at cost $\beta \geq 0$, or decide to refrain from intervening. The objective is to identify an optimal trade-off between minimizing the long-run expected avalanche size and the long-run average intervention cost. In order to analyze this optimization problem, we formulate it as a Markov decision process \citep{Puterman}. 

The state space of this MDP is the set of stable configurations $\Omega$. The action space is given by
\[
A = V \cup \{0\},
\]
where action $a \in V$ corresponds to removing a sand grain from vertex $a$ and action $0$ corresponds to no intervention. The set of feasible actions in configuration $\eta \in \Omega$ is defined as
\[
A(\eta) = V^+(\eta) \cup \{0\}.
\]
Given a configuration $\eta \in \Omega$ and an action $a \in A(\eta)$, let $\eta^a$ denote the post-decision configuration,  defined by
\[
\eta^a(x) = \begin{cases}
    \eta(x) - 1, &\text{if } a = x, \\
    \eta(x), &\text{otherwise},
\end{cases}
\qquad x \in V.
\]

Thus, $\eta^a$ is the configuration obtained immediately after applying action $a$, but before the addition of the new sand grain. The transition kernel $P: \Omega \times A \times \Omega \rightarrow [0,1]$ of the MDP is given by
\[
P(\eta' \mid \eta, a) = p(\eta^a, \eta'), \qquad \eta, \eta' \in \Omega, \quad a \in A(\eta),
\]
where $p:\Omega \times \Omega \rightarrow [0,1]$ is the transition kernel of the uncontrolled sandpile Markov chain defined in (\ref{trans_probs_uncontrolled}). 
The cost function $c: \Omega \times A \rightarrow \mathbb{R}$ consists of two components: the expected size of the next avalanche and the intervention cost. It is defined as
\begin{equation}\label{cost_function}
c(\eta, a) = \mathbb{E}[N(\eta^a)] + \beta \mathbbm{1}[a \neq 0], \qquad \eta \in \Omega, \quad a \in A(\eta).  
\end{equation}
We assume that the decision maker follows a stationary deterministic policy, which is determined by a decision rule $d: \Omega \rightarrow A$ satisfying $d(\eta) \in A(\eta)$ for every $\eta \in \Omega$. Such a policy is denoted by $\pi = d^{\infty}$. Let $\Pi$ denote the set of all stationary deterministic policies. Associated with a decision rule $d$, let $P_d: \Omega \times \Omega \rightarrow [0,1]$ and $c_d: \Omega \rightarrow \mathbb{R}$ denote the transition kernel and cost function induced by $d$, namely,
\[
P_d(\eta, \eta') = P(\eta' \mid \eta, d(\eta)), \quad \text{and } \quad c_d(\eta) = c(\eta, d(\eta)),
\]
for $\eta, \eta' \in \Omega$. We refer to the Markov chain induced by $\pi$ as the $\pi$-dynamics. 

We focus on a particular class of policies that remove sand grains from a fixed set of vertices whenever possible. We call such policies \textit{local intervention policies}. 

\begin{definition}\label{local_intervention}
    Given a graph $\mathcal{G} = (V \cup \{s\}, E)$ and a subset $W \subseteq V$, a policy $\pi_W = (d_W)^{\infty}$ is called a local intervention policy with intervention set $W$ if
    \[
    d_W(\eta) \in \begin{cases} 
        W \cap V^+(\eta) &\text{if } W \cap V^+(\eta) \neq \emptyset, \\
        \{0\}, &\text{otherwise,}
        \end{cases} \qquad \eta \in \Omega.
    \]
\end{definition}
Note that Definition \ref{local_intervention} does not specify which occupied vertex in the set $W$ is selected when there are multiple options. Hence, there may be different local intervention policies with the same intervention set. Whenever we explicitly specify the selection rule, we indicate this with a suitable superscript. 

\subsection{The average-cost optimality equations}
\label{ssec:optimality-eqns}
To assess the performance of a policy $\pi \in \Pi$, we use the average-cost criterion. We consider repeated interventions in a system that is continually driven by the addition of new sand grains, without a terminal horizon or priority given to reducing earlier costs. Our objective is therefore to achieve a sustained reduction in expected avalanche costs relative to the intervention costs. The average-cost formulation also connects the control problem to the recurrent classes and stationary behavior of the perturbed sandpile Markov chain, which are of independent interest in the theory of the Abelian sandpile model. The average-cost function, or gain, of a general stationary deterministic policy $\pi = d^{\infty} \in \Pi$ is defined as
\begin{align}
    \label{gpi_eta} g^{\pi}(\eta) &= \lim_{T \rightarrow \infty} \dfrac{1}{T} \mathbb{E}^{\pi}_{\eta}\left[\sum\limits_{t=1}^T c_d(\eta_t)\right], \quad \eta \in \Omega,
\end{align}
where $(\eta_t)_{t \geq 1}$ denotes the sandpile configuration process under the policy $\pi$. Its law and initial configuration $\eta_1 = \eta$ are indicated by the superscript and subscript in the expectation. The sandpile MDP is multichain, meaning that the Markov chain induced by a stationary deterministic policy may have multiple recurrent classes. This occurs, for example, for the always-intervene policy on the star graph, see Lemma \ref{value_always_intervene_star}. Therefore, the gain may in general depend on the initial configuration.

Given a stationary deterministic policy $\pi \in \Pi$, let $R_1, \ldots, R_k$ denote the recurrent classes of the Markov chain induced by this policy. Let $\mu_{\pi, i}$ denote the stationary distribution on $R_i$. Moreover, for an initial sandpile configuration $\eta \in \Omega$, let $\delta_i(\eta)$ denote the probability that the chain eventually enters the recurrent class $R_i$. The gain can be written as \citep{Puterman}
\begin{equation}\label{av_cost_def_2}
    g^{\pi}(\eta) = \sum\limits_{i = 1}^k \delta_i(\eta) \sum\limits_{\eta' \in R_i} \mu_{\pi, i}(\eta')c_d(\eta').
\end{equation}
Let the bias of a policy $\pi = d^{\infty} \in \Pi$ be defined as\footnote{
The definition assumes that under policy $\pi$, the controlled Markov chain is aperiodic; otherwise the limit is replaced by a Cesaro-limit (see \cite[Section 8.2.1]{Puterman}). However, we will be able to by-pass periodicity issues.
} 
\begin{equation}\label{bias}
    h^{\pi}(\eta) = \lim_{T \rightarrow \infty} \mathbb{E}^{\pi}_{\eta}\left[\sum\limits_{t=1}^{T}[c_d(\eta_t) - g^{\pi}(\eta_t)]\right].
\end{equation}
The bias can be interpreted as the expected cumulative deviation of the one-stage costs from the long-run average cost under policy $\pi$, starting from configuration $\eta$. Since the configuration space is finite, the gain $g^{\pi}$ and bias $h^{\pi}$ of a policy $\pi = d^{\infty} \in \Pi$ satisfy the following policy evaluation equations \cite[Thm.\ 8.2.6]{Puterman}
\begin{subequations}\label{policy_evaluations_equations}
\begin{align}
    \label{pol_ev_1}(I-P_d)g^{\pi} &= 0, \\
    \label{pol_ev_2}c_d - g^{\pi} + (P_d-I)h^{\pi}&= 0,
\end{align}
\end{subequations}
where $g^{\pi}$, $h^{\pi}$, and $c_d$ are regarded as column vectors indexed by the configuration space.
The policy evaluation equations (\ref{pol_ev_1}) and (\ref{pol_ev_2}) uniquely determine the gain of policy $\pi$, whereas they determine the bias only up to addition of an element of the null space of the matrix $I-P_d$. We call any function $h^{\pi}$ satisfying the equations a bias function associated with the policy $\pi$.  

A policy $\pi^*$ is called average-cost optimal if it satisfies
\begin{equation*}
    g^{\pi^*}(\eta) \leq g^{\pi}(\eta), \qquad \text{ for all } \eta \in \Omega, \quad \pi \in \Pi. 
\end{equation*}
The optimal gain is
\begin{equation*}
    g^*(\eta) = \min_{\pi \in \Pi}g^{\pi}(\eta).
\end{equation*} 
Since the state and action spaces are finite and the cost function is bounded, there exists an optimal stationary deterministic policy \cite[Thms. 9.1.4 and 9.1.8]{Puterman}. This justifies restricting our attention to this class of policies.

Let $\mathcal{F}(\Omega)$ denote the space of real-valued functions on $\Omega$. The optimal gain can be characterized by means of the following system of optimality equations \cite[Thm.\ 9.1.3]{Puterman}:
\begin{subequations}\label{optimality_equations}
\begin{align}
\label{optimality_equation_1}&\min_{a \in A(\eta)} \left\{\sum\limits_{\eta' \in \Omega} P(\eta'|\eta, a)g(\eta') - g(\eta)\right\} = 0,\\
\label{optimality_equation_2} &\min_{a \in B_{\eta}(g)} \left\{c(\eta, a) - g(\eta) + \sum\limits_{\eta' \in \Omega} P(\eta'|\eta, a)h(\eta') - h(\eta) \right\} = 0,
\end{align}
\end{subequations}
where 
\[
B_{\eta}(g) = \left\{a \in A(\eta): \sum\limits_{\eta' \in \Omega} P(\eta'|\eta, a)g(\eta') - g(\eta) = 0 \right\}.
\]
The existence of a solution to these equations is established in \cite[Thm.\ 9.1.4]{Puterman}. The function $g$ in every solution equals the optimal gain $g^*$ \cite[Thm. 9.1.3]{Puterman}. 
Conversely, if $(g, h) \in \mathcal{F}(\Omega) \times \mathcal{F}(\Omega)$ is a pair of functions that satisfies the optimality equations (\ref{optimality_equation_1}, \ref{optimality_equation_2}) and a decision rule $d^* \in D$ satisfies 
\begin{align*}
&P_{d^*}g = g, \\
&d^* \in \argmin_{d \in E(g)}\{c_d + P_d h\},
\end{align*}
where $E(g) = \{d \in D: d(\eta) \in B_{\eta}(g)\}$, then the policy $\pi^* = (d^*)^{\infty}$ is average-cost optimal \cite[Thm.\ 9.1.7]{Puterman}.

\section{Local intervention policies on general graphs}
\label{General_results_for_local_intervention_policies}

In this section, we study the controlled sandpile Markov chain under local intervention policies on general finite graphs. First, we analyze the local intervention policy for which the intervention set is empty, i.e., the policy that never intervenes. These results form the basis for the analysis of local intervention policies with general intervention sets. We provide structural results concerning the recurrent classes and stationary distributions of the Markov chains induced by local interventions. Then, we introduce the concept of topple-freeness and show that under this condition many of these results can be sharpened. Finally, under the topple-freeness condition, we provide explicit expressions for the gain and bias of local intervention policies and bounds on the immediate intervention cost $\beta$ that guarantee their optimality.

\subsection{The never-intervene policy}
\label{The_never-intervene_policy}

In this subsection, we analyze the local intervention policy $\pi_{\emptyset} = (d_{\emptyset})^{\infty}$ with empty intervention set, which we call the \textit{never-intervene policy}. Under this policy, the dynamics coincide with those of the uncontrolled sandpile Markov chain. We derive explicit expressions for its gain and bias and determine the exact threshold above which it is optimal. 

Given a graph $\mathcal{G}= (V \cup \{s\}, E)$, let 
\[
G = \Delta^{-1}
\]
denote the Green's function of $\mathcal{G}$. For any vertex $x \in V$, we define
\begin{equation}\label{def_phi}
    \phi(x) = \sum\limits_{y \in V} G(x,y).
\end{equation}
The function $\phi$ admits the following probabilistic interpretation. Let 
\[
\tau_x = \inf\{t \geq 0: X_t = x\}, \qquad x \in V \cup \{s\},
\]
denote the first hitting time of vertex $x$ of a continuous-time random walk $(X_t)_{t \geq 0}$ with jump rate 1 along each incident edge.
This walk is the continuous-time Markov chain with generator $-\Delta$, killed when it reaches $s$. Then (see e.g.~\cite[Section 4.4, equation 4.21]{Durrett_Stoch_Proc})

\begin{equation}\label{prob_interpret}
    \phi(x) = \mathbb{E}_x[\tau_s].
\end{equation}
Theorem \ref{gain_bias_never_intervene} expresses the gain and bias of the never-intervene policy in terms of the function $\phi$. 

\begin{thm}\label{gain_bias_never_intervene}
    For the policy $\pi_{\emptyset} = (d_{\emptyset})^{\infty}$, where $d_{\emptyset}(\eta) = 0$ for all $\eta \in \Omega$, the gain is given by
    \begin{equation}\label{g_pi0} 
        g^{\pi_{\emptyset}}(\eta) = \dfrac{1}{|V|}\sum_{x \in V} \phi(x), \qquad \eta \in \Omega.
    \end{equation}
    Moreover, there exists a constant $C \in \mathbb{R}$ such that
    \begin{equation}\label{h_pi0}
        h^{\pi_{\emptyset}}(\eta) = \sum\limits_{x \in V} (\eta(x) - \deg(x) + 1)\phi(x) + C, \qquad \eta \in \Omega.
    \end{equation} 
\end{thm}
\begin{proof}
    Let $g^{\pi_{\emptyset}}$ and $h^{\pi_{\emptyset}}$ be given by (\ref{g_pi0}) and (\ref{h_pi0}). Since $g^{\pi_{\emptyset}}$ is constant, (\ref{pol_ev_1}) holds immediately. It remains to verify that $g^{\pi_{\emptyset}}$ and $h^{\pi_{\emptyset}}$ satisfy (\ref{pol_ev_2}). From the definition of the stabilization process, we have
    \begin{equation*}
        (\alpha_x \eta)^{\circ}(y) = (\alpha_x \eta)(y) - \sum\limits_{z \in V} n_z(\eta, x)\Delta_{zy}, \qquad y \in V.
    \end{equation*}
    Hence, we obtain
    \begin{align*}
        h^{\pi_{\emptyset}}((\alpha_x \eta)^{\circ}) - h^{\pi_{\emptyset}}(\eta) &= \sum\limits_{y \in V} \left[\mathbbm{1}_x(y) - \sum\limits_{z \in V} n_z(\eta, x)\Delta_{zy}\right] \phi(y) \\
        &= \phi(x) - \sum\limits_{z \in V}n_z(\eta, x) \sum\limits_{y \in V} \Delta_{zy}\phi(y). 
    \end{align*}
    Since $G = \Delta^{-1}$ and $\phi = G\mathbf{1}$, we have
    \begin{equation*}
        \sum\limits_{y \in V}\Delta_{zy}\phi(y) = 1, \qquad \text{for all } z \in V.
    \end{equation*}
    Hence,
    \begin{equation*}
        h^{\pi_{\emptyset}}((\alpha_x \eta)^{\circ}) - h^{\pi_{\emptyset}}(\eta) = \phi(x) - n(\eta, x).
    \end{equation*}
    Thus, we obtain
    \begin{equation*}
        ((P_{d_{\emptyset}} - I)h^{\pi_{\emptyset}})(\eta) = \dfrac{1}{|V|} \sum\limits_{x \in V} \phi(x) - \mathbb{E}[N(\eta)].
    \end{equation*}
    Using expression (\ref{g_pi0}) and the fact that 
    \begin{equation*}
        c_{d_{\emptyset}}(\eta) = \mathbb{E}[N(\eta)] \quad \text{for all } \eta \in \Omega,
    \end{equation*}
    we conclude that 
    \begin{equation*}
        c_{d_{\emptyset}} - g^{\pi_{\emptyset}} + (P_{d_{\emptyset}} - I)h^{\pi_{\emptyset}} = 0.
    \end{equation*}
    Since the Markov chain induced by $\pi_{\emptyset}$ has a unique recurrent class, the policy evaluation equations determine the bias up to an additive constant. Hence, every bias function associated with $\pi_{\emptyset}$ has the form (\ref{h_pi0}).
\end{proof}

\begin{remark}
   Formula (\ref{g_pi0}) can alternatively be obtained directly from the definition of the gain (\ref{gpi_eta}) and Dhar's formula~\citep{Dhar}. 
\end{remark}

Using the expressions (\ref{g_pi0}) and (\ref{h_pi0}) for the gain and bias of the never-intervene policy, we characterize the values of the immediate intervention cost $\beta$ for which this policy is optimal.

\begin{thm}\label{nothing_optimal}
    Consider the policy $\pi_{\emptyset} = (d_{\emptyset})^{\infty}$, where $d_{\emptyset}(\eta) = 0$ for all $\eta \in \Omega$. Policy $\pi_{\emptyset}$ is optimal if and only if 
    \begin{equation}\label{cond_nothing_optimal}
        \beta \geq \max\{\phi(a) : a \in V,\ \deg(a)\geq 2\}.
    \end{equation}
\end{thm}
\begin{proof}
    In order to prove sufficiency, we verify that $(g^{\pi_{\emptyset}}, h^{\pi_{\emptyset}})$ as defined in expressions (\ref{g_pi0}) and (\ref{h_pi0}) satisfy the average-cost optimality equations (\ref{optimality_equation_1}) and (\ref{optimality_equation_2}) whenever condition (\ref{cond_nothing_optimal}) holds. The validity of equation (\ref{optimality_equation_1}) follows immediately from the fact that $g^{\pi_{\emptyset}}$ is constant over all sandpile configurations. We continue to prove equation (\ref{optimality_equation_2}). Since $g^{\pi_{\emptyset}}$ is constant, the set $B_{\eta}(g^{\pi_{\emptyset}})$ equals $A(\eta)$ for every $\eta \in \Omega$. Observe that the expressions (\ref{g_pi0}) and (\ref{h_pi0}) satisfy the policy evaluation equation (\ref{pol_ev_2}) for policy $\pi_{\emptyset}$. Hence, action $0$ attains the minimum in equation (\ref{optimality_equation_2}) if
    \begin{equation}\label{optimality_equation_exp}
        c(\eta, 0) + \sum\limits_{\eta' \in \Omega} P(\eta'|\eta, 0)h^{\pi_{\emptyset}}(\eta') \leq c(\eta, a) + \sum\limits_{\eta' \in \Omega}P(\eta'|\eta, a)h^{\pi_{\emptyset}}(\eta'), \qquad \text{for all } \eta \in \Omega, \quad a \in A(\eta). 
    \end{equation}
    Using the facts that
    \begin{equation*}
        c(\eta, a) = \beta \mathbbm{1}\{a \neq 0\} + c(\eta^a, 0),
    \end{equation*}
    and 
    \[
    P(\eta'|\eta, a) = P(\eta'|\eta^a, 0),
    \]
    together with the policy evaluation equation (\ref{pol_ev_2}) for $\pi_{\emptyset}$, expression (\ref{optimality_equation_exp}) is equivalent to
    \[
    g^{\pi_{\emptyset}}+h^{\pi_{\emptyset}}(\eta) \leq \beta \mathbbm{1}\{a \neq 0\} + g^{\pi_{\emptyset}}+h^{\pi_{\emptyset}}(\eta^a), \qquad \eta \in \Omega, \quad a \in A(\eta).
    \]
    It follows that expression (\ref{optimality_equation_exp}) holds if and only if
    \[
    \beta \mathbbm{1}\{a \neq 0\} \geq h^{\pi_{\emptyset}}(\eta) - h^{\pi_{\emptyset}}(\eta^a), \qquad \eta \in \Omega, \quad a \in A(\eta).
    \]
    Since $\eta$ and $\eta^a$ differ only at site $a$, expression (\ref{h_pi0}) ensures that the optimality equations hold whenever
    \begin{equation}
        \beta \geq \phi(a) \quad \text{for all } a \in V \text{ for which there exists } \eta \in \Omega \text{ with } \eta(a) > 0. 
    \end{equation}
    Since there exists a stable configuration $\eta \in \Omega$ with $\eta(a) > 0$ if and only if $\deg(a) \geq 2$, the above condition is equivalent to (\ref{cond_nothing_optimal}), which establishes sufficiency.

    We proceed to prove necessity. Given $\eta \in \Omega$, $a \in A(\eta)$, and $x \in V$, the definition of stabilization yields
    \[
    (\alpha_x \eta^a)^{\circ}(y) = \eta^a(y) + \mathbbm{1}_x(y) - \sum\limits_{u \in V}n_u(\eta^a, x)\Delta_{uy}.
    \]
    Using the explicit expression (\ref{h_pi0}) and the identity $(\Delta \phi)(u) = 1$, we obtain
    \[
    h^{\pi_{\emptyset}}((\alpha_x \eta^a)^{\circ}) - h^{\pi_{\emptyset}}(\eta) = \phi(x) - n(\eta^a, x) - \phi(a)\mathbbm{1}\{a \neq 0\},
    \]
    where we set $\phi(0) = 0$. Therefore,
    \[
    n(\eta^a, x) + \beta \mathbbm{1}\{a \neq 0\} = \phi(x) - \left[h^{\pi_{\emptyset}}((\alpha_x\eta^a)^{\circ}) - h^{\pi_{\emptyset}}(\eta)\right] + (\beta - \phi(a))\mathbbm{1}\{a \neq 0\}.
    \]
    Taking the expectation over the addition site now yields
    \begin{equation}\label{cost_exp}
    c(\eta, a) = g^{\pi_{\emptyset}} + h^{\pi_{\emptyset}}(\eta) - \sum\limits_{\eta' \in \Omega}P(\eta'|\eta, a)h^{\pi_{\emptyset}}(\eta') + (\beta - \phi(a))\mathbbm{1}\{a \neq 0\}.
    \end{equation}
    Suppose that $\beta < \phi(a)$ for some $a \in V$ with $\deg(a) \geq 2$. Consider the stationary deterministic policy $\pi_a = (d_a)^{\infty}$, where 
    \begin{equation}\label{pi_a-policy_def}
    d_a(\eta) = \begin{cases}
        a, &\text{if } \eta(a) > 0, \\
        0, &\text{if } \eta(a) = 0.
    \end{cases}
    \end{equation}
    Let $R$ denote a recurrent class of the Markov chain induced by $\pi_a$. We claim that there is at least one configuration $\eta \in R$ with $\eta(a) > 0$. Indeed, suppose that $\eta(a) = 0$ for all $\eta \in R$. Then, for any $\eta \in R$, there is no intervention by the policy. If the next grain is added at $a$, which occurs with positive probability, we obtain a configuration $\eta'$ with $\eta'(a) = 1$, since $\deg(a) \geq 2$. By assumption, $\eta' \notin R$, which contradicts the closedness of $R$. Let $\mu_{\pi_a, R}$ denote the stationary distribution of the Markov chain induced by $\pi_a$ on $R$. Since $R$ is irreducible and contains at least one configuration $\zeta$ satisfying $\zeta(a) > 0$, we have
    \[
    \sum_{\zeta \in R}\mu_{\pi_a, R}(\zeta)\mathbbm{1}\{\zeta(a) > 0\} > 0.
    \]
    Taking expectations in (\ref{cost_exp}) with respect to $\mu_{\pi_a, R}$ and using its invariance under the $\pi_a$-dynamics, we obtain
    \[
    g^{\pi_a}(\eta) = g^{\pi_{\emptyset}} + \sum\limits_{\zeta \in R}\mu_{\pi_a, R}(\zeta)\mathbbm{1}\{\zeta(a) > 0\}(\beta - \phi(a)) < g^{\pi_{\emptyset}}(\eta), \qquad \eta \in R.
    \]
    Thus, $\pi_{\emptyset}$ is not optimal whenever $\beta < \phi(a)$ for some $a \in V$ with $\deg(a) \geq 2$. This establishes necessity.
\end{proof}

\subsection{Structural results for local intervention policies}
\label{Structural_results_for_local_intervention_policies}

In this subsection, we show that, although a local intervention policy perturbs the sandpile dynamics, part of the structure of the uncontrolled sandpile Markov chain is preserved. In particular, we relate the recurrent configurations of the controlled Markov chain to those of the uncontrolled sandpile dynamics on the graph obtained by identifying all vertices in $W$ with the sink. For singleton intervention sets, we also establish a relation between the stationary distributions. 

Given a graph $\mathcal{G} = (V \cup \{s\}, E)$, let
\[
\tau_W = \inf\{t \geq 0: X_t \in W\}, \qquad W \subseteq V \cup \{s\},
\]
denote the first hitting time of the set $W$ by a simple random walk $(X_t)_{t \geq 0}$. For a singleton set $\{x\}$, we simply write $\tau_x$  and we set $\tau_{\emptyset} = \infty$. Furthermore, for $W \subseteq V$ and $x \in V \setminus W$, let 
\[
p_W(x) = \mathbb{P}_x(\tau_W < \tau_s)
\]
denote the probability that a simple random walk started at $x$ hits $W$ before hitting the sink. 

Let $\widetilde{\mathcal{G}}_W$ denote the graph obtained from $\mathcal{G}$ by contracting all vertices in $W \cup \{s\}$ into a single sink vertex. Every edge from a vertex in $V \setminus W$ to a vertex in $W \cup \{s\}$ becomes an edge to the contracted sink in this procedure. Edges with both endpoints in $W \cup \{s\}$ are discarded. We refer to this graph as the contracted graph. Let $\widetilde{\Omega}_W$ denote the set of stable configurations on $\widetilde{\mathcal{G}}_W$, and let $\widetilde{R}_W$ denote the unique recurrent class of the uncontrolled sandpile Markov chain on $\widetilde{\mathcal{G}}_W$.  We let $\phi_W$ denote the analogue of the function $\phi$ defined in expression (\ref{def_phi}) for the contracted graph $\widetilde{\mathcal{G}}_W$. Furthermore, let $R_{\pi_W}$ denote the set of recurrent configurations of the Markov chain on $\mathcal{G}$ induced by the policy $\pi_W$, which may be the union of several recurrent classes. Given a stationary deterministic policy $\pi \in \Pi$ for the sandpile MDP on the graph $\widetilde{\mathcal{G}}_W$, we let $\widetilde{g}^{\pi}_W$ and $\widetilde{h}^{\pi}_W$ denote its gain and bias functions, respectively. For $\eta \in \Omega$, we write
\[
\xi = \eta|_{V \setminus W} \in \widetilde{\Omega}_W
\]
for the restriction of $\eta$ to the graph $\widetilde{\mathcal{G}}_W$. 

We first establish a relation between the recurrent configurations under a local intervention policy $\pi_W$ and those of the uncontrolled sandpile model on the contracted graph $\widetilde{\mathcal{G}}_{W}$ in Theorem \ref{rec_conf_pi_W}. 

\begin{thm}\label{rec_conf_pi_W}
Given a graph $\mathcal{G} = (V \cup \{s\}, E)$ and a subset $W \subseteq V$, let $\pi_W = (d_W)^{\infty}$ denote a local intervention policy. We have
\begin{equation*}
    \{\eta|_{V \setminus W}: \eta \in R_{\pi_W}\} \subseteq \widetilde{R}_W.
\end{equation*}
\end{thm}
\begin{proof}
    Consider a configuration $\eta \in R_{\pi_W}$ and let $\xi = \eta|_{V \setminus W}$. Let $\eta'$ denote the result of one transition under the policy $\pi_W$, obtained from an intervention by the policy, addition at a vertex $x \in V$ and subsequent stabilization, i.e. $\eta' = (\alpha_x \eta^{d_W(\eta)})^{\circ}$. Recall that $n_y( \eta^{d_W(\eta)}, x)$ denotes the number of times a vertex $y \in V$ topples during this transition.  We define
    \[
    \gamma_u = \sum\limits_{y \in W}n_y(\eta^{d_W(\eta)}, x)b_{yu}, \qquad u \in V \setminus W,
    \]
    i.e., the total number of grains sent from $W$ to $u$ during the transition.
    By the Abelian property, the restriction $\xi' := \eta'|_{V \setminus W}$ can be obtained from $\xi$ by first adding a grain at $x$ if $x \notin W$, then adding $\gamma_u$ grains at each vertex $u \in V\setminus W$ and finally stabilizing the resulting configuration on the contracted graph $\widetilde{\mathcal{G}}_W$. Letting $\alpha_u^k$ denote the operation of adding $k$ grains to vertex $u$, for $k = 1, 2, \ldots,$ and $\alpha_u^0$ the identity operator, the configuration $\xi'$ can be written as 
    \begin{equation*}
        \xi' = \left(\left(\prod_{u \in V \setminus W} \alpha_u^{\gamma_u}\right) \alpha_x^{\mathbbm{1}\{x \notin W\}} \xi\right)^{\circ}.
    \end{equation*}
    
    Hence, the transition from $\xi$ to $\xi'$ is realized by a finite sequence of additions and relaxations on $\widetilde{\mathcal{G}}_W$. Since $\widetilde{R}_W$ is closed under addition and relaxation operators on $\widetilde{\mathcal{G}}_W$, it follows that 
    \begin{equation}\label{closed}
        \xi \in \widetilde{R}_W \quad \Rightarrow \quad \xi' \in \widetilde{R}_W.
    \end{equation}
    Now, suppose that the restriction $\xi \notin \widetilde{R}_W$. Since $\widetilde{R}_W$ is the unique recurrent class of the uncontrolled sandpile dynamics on $\widetilde{\mathcal{G}}_W$, there exists a finite sequence $x_1, \ldots, x_k \in V\setminus W$ such that,
    \[
        (\alpha_{x_k}\cdots \alpha_{x_1} \xi)^{\circ} \in \widetilde{R}_W.
    \]
    Consider the process on $\mathcal{G}$ under the $\pi_W$-dynamics. The same sequence of additions at $x_1, \ldots, x_k \in V\setminus W$ may cause topplings of vertices in the set $W$. For $y \in W$, let $N_y$ denote the total number of topplings of $y$ during these transitions. By the Abelian property, the resulting configuration restricted to $V\setminus W$ can be  obtained by first performing the additions at $x_1, \ldots, x_k$ and stabilizing on $\widetilde{\mathcal{G}}_W$, then adding $\sum_{y \in W}N_y b_{yu}$ grains to each vertex $u \in V \setminus W$, and stabilizing again. Hence, it can be written as 
    \[
    \left(\left(\prod_{u \in V\setminus W} \alpha_u^{\sum\limits_{y \in W}N_yb_{yu}}\right)(\alpha_{x_k}\cdots \alpha_{x_1} \xi)^{\circ}\right)^{\circ}.
    \]
    Since $(\alpha_{x_k} \cdots \alpha_{x_1}\xi)^{\circ} \in \widetilde{R}_W$ and $\widetilde{R}_W$ is closed under addition and stabilization, this configuration also belongs to $\widetilde{R}_W$. Hence, we have shown that the set
    \[
        \{\zeta \in \Omega: \zeta|_{V \setminus W} \in \widetilde{R}_W\}
    \]
    is closed under the $\pi_W$-dynamics by (\ref{closed}) and that it can be reached from $\eta$ with positive probability. Since $\xi \notin \widetilde{R}_W$, this set does not contain $\eta$. Therefore, there is a positive probability that the process leaves $\eta$ and can never return, contradicting the recurrence of $\eta$. Thus, if $\eta \in R_{\pi_W}$, then $\eta|_{V \setminus W} \in \widetilde{R}_W$. 
\end{proof}
The reverse inclusion in Theorem \ref{rec_conf_pi_W} does not hold in general, as the following counterexample demonstrates.

\paragraph{Counterexample to the reverse inclusion.}
Consider a graph $\mathcal{G} = (V \cup \{s\}, E)$, where $V = \{x, y\}$ and $E = \{(x, y), (y, s)\}$. Then, the process induced by the local intervention policy $\pi_x$ has a unique recurrent class $R_{\pi_x} = \{(0,1)\}$.
Consequently,
\begin{equation*}
    \{\eta|_{V \setminus \{x\}}: \eta \in R_{\pi_x}\} = \{1\}.
\end{equation*}
On the other hand, the uncontrolled sandpile dynamics on the contracted graph $\widetilde{\mathcal{G}}_x$ has recurrent class $\widetilde{R}_x = \{0,1\}$.
This implies
\begin{equation*}
    \{\eta|_{V \setminus \{x\}}: \eta \in R_{\pi_x}\} \subsetneq \widetilde{R}_x.
\end{equation*}
Hence, the reverse inclusion in Theorem \ref{rec_conf_pi_W} does not hold in general. 
\qed
\\
\\
We proceed to investigate the stationary distributions of the Markov chains induced by local intervention policies with singleton intervention sets and relate them to those of the uncontrolled sandpile dynamics on the corresponding contracted graph.    
Let $\widetilde{\mu}_x$ denote the stationary distribution of the uncontrolled sandpile dynamics on the contracted graph $\widetilde{\mathcal{G}}_x$. A toppling of $x$ in the original graph sends $b_{xy}$ grains to each vertex $y \in V \setminus \{x\}$. Therefore, the effect of such a toppling on the restricted configuration, followed by stabilization on the contracted graph, is represented by the operator $B_x: \widetilde{\Omega}_x \rightarrow \widetilde{\Omega}_x$ defined as
\begin{equation*}
    B_x \xi = \left(\left[\prod\limits_{y \in V \setminus \{x\}} \alpha_y^{b_{xy}}\right]\xi \right)^{\circ}.
\end{equation*}
Observe that the restriction of $B_x$ to $\widetilde{R}_x$ is an element of the sandpile group $\mathcal{S}(\widetilde{\mathcal{G}}_x)$. Consequently, it generates a cyclic subgroup $\langle B_x \rangle \leq \mathcal{S}(\widetilde{\mathcal{G}}_x)$. We denote by $C_1, \ldots, C_N$ the orbits of $B_x$ on $\widetilde{R}_x$.  Since recurrent configurations can be identified with sandpile group elements,  these orbits correspond to the cosets of the cyclic subgroup $\langle B_x \rangle$ in $\mathcal{S}(\widetilde{\mathcal{G}}_x)$. The following result shows that the projection of any stationary distribution $\mu_{\pi_x}$ onto $V \setminus \{x\}$ assigns equal mass to these orbits. 

\begin{thm}\label{uniform_cosets}
    Given a graph $\mathcal{G} = (V \cup \{s\}, E)$ and a vertex $x \in V$, let $\mu_{\pi_x}$ denote any stationary distribution of the $\pi_x$-dynamics and let $\bar{\mu}_{\pi_x}$ denote the projection of $\mu_{\pi_x}$ onto $\widetilde{\Omega}_x$. Then, $\bar{\mu}_{\pi_x}$ assigns equal probability to the orbits $C_1, \ldots, C_N$ of $B_x$ on $\widetilde{R}_x$, i.e.,
    \begin{equation*}
        \bar{\mu}_{\pi_x}(C_i) = \dfrac{1}{N}, \quad \text{for all } i = 1, \ldots, N.
    \end{equation*}
\end{thm}
\begin{proof}
Let $\eta_0 \sim \mu_{\pi_x}$, and let $\xi_0 = \eta_0|_{\mathcal{G} \setminus \{x\}}$ denote its restriction to $\mathcal{G} \setminus \{x\}$. Observe that Theorem \ref{rec_conf_pi_W} implies that $\xi_0 \in \widetilde{\mathcal{R}}_x$. Consider the following coupling of the $\pi_x$-dynamics on $\mathcal{G}$ and the uncontrolled dynamics on the contracted graph $\widetilde{\mathcal{G}}_x$: at each time step, we select a vertex $y \in V$ uniformly at random. If $y \neq x$, we perform the transition under $\pi_x$ in the first process, i.e., first remove a grain from $x$ if possible, then add at $y$ and relax. In the second process, we add at $y$ and relax on $\widetilde{\mathcal{G}}_x$. If $y = x$, we perform the $\pi_x$-transition in the first process and do nothing in the second process. Let $\eta_t \in \Omega$ denote the configuration obtained after running the $\pi_x$-dynamics for $t$ steps and let $\xi_t$ denote the result of running the uncontrolled dynamics on the contracted graph $\widetilde{\mathcal{G}}_x$ under this coupling. Let $k_t$ denote the total number of topplings of $x$ up to time $t$ in the evolution under the $\pi_x$ dynamics. Then, by the Abelian property,
\begin{equation*}
    \eta_t|_{\mathcal{G} \setminus \{x\}} = B_x^{k_t}\xi_t.
\end{equation*}
This implies 
\begin{equation*}
    \eta_t|_{\mathcal{G} \setminus \{x\}} \in C_i \Leftrightarrow \xi_t \in C_i, \quad i = 1, \ldots, N,
\end{equation*}
since $C_i$ is invariant under $B_x$ for each $i = 1, \ldots, N$. We now obtain
\begin{equation*}
    \bar{\mu}_{\pi_x}(C_i) = \lim_{T \rightarrow \infty} \mathbb{P}(\eta_T|_{\mathcal{G} \setminus \{x\}} \in C_i) = \lim_{T \rightarrow \infty} \mathbb{P}(\xi_T \in C_i) = \widetilde{\mu}_x(C_i) = \dfrac{|C_i|}{|\widetilde{\mathcal{R}}_x|},
\end{equation*}
since the stationary distribution of the second process is the same as that of the uncontrolled sandpile dynamics on the contracted graph. Since all cosets of a subgroup have the same cardinality, this yields
\begin{equation*}
    \bar{\mu}_{\pi_x}(C_i) = \dfrac{1}{N}, \quad \text{for all } i = 1, \ldots, N.
\end{equation*}
\end{proof}

\subsection{Topple-free intervention sets}
\label{Topple-free_intervention_sets}

In this subsection, we introduce the notion of topple-freeness. We show that under this condition, the results in Section~\ref{Structural_results_for_local_intervention_policies} can be considerably sharpened. For local intervention policies with topple-free intervention sets, we derive explicit expressions for the gain and bias and establish sufficient conditions for optimality.

\begin{definition}\label{topple_free}
    Consider a graph $\mathcal{G} = (V \cup \{s\}, E)$ and a stationary deterministic policy $\pi$. A set $W \subseteq V$ is \textit{topple-free} under $\pi$ if, under every stationary distribution of the Markov chain induced by $\pi$, the probability that any vertex $x \in W$ topples during a transition is zero. In particular, we say that a vertex $x \in V$ is topple-free under $\pi$ if the singleton set $\{x\}$ is topple-free under $\pi$. 
\end{definition}

\subsubsection{Recurrent classes and stationary distributions under topple-freeness}
\label{Recurrent_classes_and_stationary_distributions_under_topple-freeness}

Since vertices in a topple-free intervention set $W$ do not topple during transitions from recurrent configurations of the Markov chain induced by the corresponding local intervention policy, more of the structure of the uncontrolled sandpile dynamics is preserved. We use this property to strengthen the structural results of the preceding subsection.

The following result shows that the reverse inclusion in Theorem \ref{rec_conf_pi_W} does hold under the topple-freeness condition. It shows that the projection of each recurrent class of the controlled dynamics onto the contracted graph equals the unique recurrent class of the uncontrolled dynamics on that graph.

\begin{thm}\label{rec_conf_pi_W_rev}
    Consider a graph $\mathcal{G} = (V \cup \{s\}, E)$ and a local intervention policy $\pi_W$ with intervention set $W \subseteq V$. If $W$ is topple-free under the policy $\pi_W$, then each recurrent class $R$ of the Markov chain induced by $\pi_W$ satisfies
    \begin{equation}\label{rec_eq_W}
        \{\eta|_{V \setminus W}: \eta \in R\} = \widetilde{R}_W.
    \end{equation}
\end{thm}
\begin{proof}
    Let $R$ denote any recurrent class of the Markov chain induced by $\pi_W$. Let 
    \[
        \widetilde{R} = \{\eta|_{V \setminus W}: \eta \in R\},
    \]
    i.e., the set of restrictions of configurations in $R$ to $V\setminus W$. By Theorem \ref{rec_conf_pi_W}, we have
    \[
        \widetilde{R} \subseteq \widetilde{R}_W.
    \]
    Consider a configuration $\eta \in R$ and let $\xi = \eta|_{V \setminus W}$ denote its restriction to $V\setminus W$. Since $W$ is topple-free under $\pi_W$, no vertex in $W$ topples during any transition from a configuration in $R$. Hence,
    \[
    ((\alpha_x \eta^{d_W(\eta)})^{\circ})|_{V \setminus W} = (\alpha_x \xi)^{\circ}, \qquad \text{for all } x \in V \setminus W,
    \]
    where addition and stabilization on the right-hand side are performed on the contracted graph $\widetilde{\mathcal{G}}_W$.
    
    Since $R$ is closed under the $\pi_W$-dynamics, it follows that $\widetilde{R}$ is closed under addition at every vertex $x \in V \setminus W$ and subsequent stabilization. 
    Since $\widetilde{R}$ is a nonempty subset of $\widetilde{R}_W$ and the uncontrolled sandpile model on $\widetilde{\mathcal{G}}_W$ is irreducible on $\widetilde{R}_W$, this implies
    \[
    \widetilde{R} = \widetilde{R}_W.
    \]
    This proves the equality (\ref{rec_eq_W}).    
\end{proof}

If the vertex $x$ is topple-free under $\pi_x$, then, we can draw a stronger conclusion about the stationary distribution than in Theorem \ref{uniform_cosets}. Indeed, in stationarity, no topplings of $x$ contribute additional grains to the contracted graph. Therefore, the projected process follows the uncontrolled dynamics on the contracted graph, with an additional holding step whenever the addition site is $x$. Consequently, every projected stationary distribution coincides with the stationary distribution of the uncontrolled dynamics on the contracted graph and is therefore uniform on $\widetilde{R}_x$, as stated in the following theorem. 
\begin{thm}
    Let $\mathcal{G} = (V \cup \{s\}, E)$ and let $x \in V$ be topple-free under the policy $\pi_x$. For any stationary distribution $\mu_{\pi_x}$ of the $\pi_x$-dynamics, let $\bar{\mu}_{\pi_x}$ denote its projection onto the configuration space of $\mathcal{G}\setminus \{x\}$. Then, $\bar{\mu}_{\pi_x}$ is the uniform distribution on $\widetilde{R}_x$. 
\end{thm}
\begin{proof}
    The fact that $x$ is topple-free under the policy $\pi_x$ implies that the number $k_t$ in the proof of Theorem \ref{uniform_cosets} equals $0$ for each $t$. Hence, under the coupling introduced in this proof, we have
    \[
    \eta_t|_{\mathcal{G}\setminus\{x\}} = \xi_t, \qquad \text{for all } t \geq 0.
    \]
    Therefore, 
    \[
    \bar{\mu}_{\pi_x}(\xi) = \lim_{T \rightarrow \infty} \mathbb{P}(\eta_T\mid_{\mathcal{G}\setminus \{x\}} = \xi) = \lim_{T \rightarrow \infty} \mathbb{P}(\xi_T = \xi) = \widetilde{\mu}_x(\xi) = \dfrac{1}{|\widetilde{R}_x|}
    \]
    for every $\xi \in \widetilde{R}_x$.
\end{proof}

\subsubsection{Gain, bias and optimality}
\label{Gain_bias_optimality}

We first express the gain and bias of a local intervention policy $\pi_W$ with topple-free intervention set $W$ in terms of those of the never-intervene policy on the contracted graph. Theorem \ref{general_optimal_policies} uses these expressions to establish sufficient conditions for optimality of this policy. 

\begin{thm}\label{gain_bias_pi_W_topple_free}
    Consider a graph $\mathcal{G} = (V \cup \{s\}, E)$ and a local intervention policy $\pi_W$ with intervention set $W \subseteq V$. Suppose that $W$ is topple-free under $\pi_W$. Then, for every $\eta \in \Omega$, the gain is given by
    \begin{equation}\label{local_gain}
        g^{\pi_W}(\eta) = \dfrac{|V|-|W|}{|V|}\widetilde{g}^{\pi_{\emptyset}}_W + \dfrac{1}{|V|}\left(|W| + \sum\limits_{x \in V\setminus W} p_W(x)\right)\beta,
    \end{equation}
    where $\widetilde{g}^{\pi_{\emptyset}}_W$ denotes the gain of the never-intervene policy on $\widetilde{\mathcal{G}}_W$.
    Moreover, for each recurrent class $R$ of the Markov chain induced by $\pi_W$, there exists a constant $C_R \in \mathbb{R}$ such that, for $\eta \in R$ and $\xi = \eta|_{V \setminus W}$, 
    \begin{equation}\label{local_bias}
        h^{\pi_W}(\eta) = \widetilde{h}^{\pi_{\emptyset}}_W(\xi) + \left[\sum\limits_{x \in W}\eta(x) + \sum_{x \in V \setminus W} p_W(x)\xi(x)\right]\beta + C_R,
    \end{equation}
    where $\widetilde{h}^{\pi_{\emptyset}}_W$ denotes the bias of the never-intervene policy on $\widetilde{\mathcal{G}}_W$.
\end{thm}

In the proof of this theorem, we are going to use a coupled random walks representation of the sandpile dynamics that we now describe (this is implicit in \cite[Section 4]{Dhar_survey}, and reference [19] therein). Let us distinguish the grains in the pile, and whenever a toppling at $x \in V$ occurs, let us choose the destinations of the $\deg_G(x)$ redistributed grains according to a uniformly random permutation among all $(\deg_G(x))!$ possibilities. Whenever a grain moves, it takes a simple random walk step on $G$. The grains together, if observed only at the times they participate in a toppling, form a coupled system of simple random walks. In particular, the number of random walk steps a grain takes is exactly the number of topplings it was involved in, and the probability that a grain added at $x \in V \setminus W$ enters $W$ equals $p_W(x)$. This coupling is useful to us, because the contributions of individual grains enter additively in the long-run cost, and hence the correlations between the walks do not play a role.

\begin{proof}[{Proof of Theorem \ref{gain_bias_pi_W_topple_free}.}]
    From the definition of the cost function (\ref{cost_function}), it follows that $g^{\pi_W}$ and $h^{\pi_W}$ are of the form
    \begin{align}
        g^{\pi_W}(\eta) \label{local_gain_decomp}&= g^{\pi_W}_{\text{av}}(\eta) + g^{\pi_W}_{\text{int}}(\eta)\beta, \\
        h^{\pi_W}(\eta) \label{local_bias_decomp}&= h^{\pi_W}_{\text{av}}(\eta) + h^{\pi_W}_{\text{int}}(\eta)\beta,
    \end{align}
    where $g^{\pi_W}_{\text{av}}$ and $h^{\pi_W}_{\text{av}}$ arise from the avalanche costs, while $g^{\pi_W}_{\text{int}}$ and $h^{\pi_W}_{\text{int}}$ arise from the intervention costs. Let $(X^{\pi_W}_t)_{t \geq 1}$ denote the Markov chain induced by the policy $\pi_W$, started from a recurrent configuration, and let $(\widetilde{X}^W_t)_{t \geq 1}$ denote a lazy version of the uncontrolled Markov chain on the contracted graph $\widetilde{\mathcal{G}}_W$ with holding probability $|W|/|V|$ at each state. We couple the two processes using the same sequence of additions. Whenever a sand grain is added at a vertex in $W$, we perform the addition only in $X^{\pi_W}$, while leaving $\widetilde{X}^W$ unchanged. Since $W$ is topple-free under $\pi_W$ and the process is started from a recurrent configuration, no vertex in $W$ topples during a transition under $\pi_W$. Hence, under this coupling, the restriction of the controlled process to $V \setminus W$ evolves according to a lazy version of the uncontrolled sandpile dynamics on $\widetilde{\mathcal{G}}_W$, and the avalanche sizes in the two processes coincide at every time step. Since in the lazy uncontrolled process on $\widetilde{\mathcal{G}}_W$, an addition occurs with probability $(|V| - |W|)/|V|$, it follows that
    \begin{equation}\label{local_gain_av}
        g^{\pi_W}_{\text{av}}(\eta) = \dfrac{|V|-|W|}{|V|}\widetilde{g}^{\pi_{\emptyset}}_W(\xi), \qquad \eta \in R_{\pi_W}.
    \end{equation}
    Likewise, we show that for any recurrent class $R$ of the Markov chain induced by the policy $\pi_W$, there exists a constant $C_R$ such that
    \begin{equation}\label{local_bias_av}
        h^{\pi_W}_{\text{av}}(\eta) = \widetilde{h}^{\pi_{\emptyset}}_W(\xi) + C_R, \qquad \eta \in R,
    \end{equation}
    where $\xi = \eta|_{V \setminus W}$. To verify (\ref{local_bias_av}), observe that equation (\ref{pol_ev_2}) for the never-intervene policy on $\widetilde{\mathcal{G}}_W$ on the set of recurrent configurations reads
    \begin{equation*}
        \left[\mathbb{E}[N(\xi)] - \widetilde{g}^{\pi_{\emptyset}}_W + \sum\limits_{\xi' \in \widetilde{R}_W}P_{d_{\emptyset}}(\xi, \xi')\widetilde{h}_W^{\pi_{\emptyset}}(\xi') - \widetilde{h}_W^{\pi_{\emptyset}}(\xi)\right] = 0, \qquad \xi \in \widetilde{R}_W.
    \end{equation*}
    Let $R$ denote a recurrent class of the Markov chain induced by $\pi_W$ and let $\eta \in R$, with $\xi = \eta|_{V \setminus W}$. By Theorem \ref{rec_conf_pi_W}, we have $\xi \in \widetilde{R}_W$. Since $W$ is topple-free under $\pi_W$, we have
    \begin{align}
        \mathbb{E}[N( \eta^{d_W(\eta)})] \label{rel_1}&= \dfrac{|V|-|W|}{|V|} \mathbb{E}[N(\xi)], \\
        \sum_{\eta'\,:\, \eta'|_{V\setminus W}=\xi'}
        P_{d_W}(\eta,\eta') \label{rel_2}&= \frac{|V|-|W|}{|V|}P_{\widetilde{d}}(\xi,\xi') + \frac{|W|}{|V|}\mathbf{1}_{\{\xi'=\xi\}}, \qquad \xi' \in \widetilde{R}_W.
    \end{align}
    Combining relations (\ref{rel_1}) and (\ref{rel_2}) with (\ref{local_gain_av}), we obtain
    \begin{align*}
        &\mathbb{E}[N(\eta^{d_W(\eta)})] - g^{\pi_W}_{\text{av}}(\eta) + \sum\limits_{\eta' \in R}P_{d_W}(\eta, \eta')h^{\pi_W}_{\text{av}}(\eta') - h^{\pi_W}_{\text{av}}(\eta) \\
        &= \dfrac{|V|-|W|}{|V|}\left[\mathbb{E}[N(\xi)]-\widetilde{g}_W^{\pi_{\emptyset}} + \sum\limits_{\xi' \in \widetilde{R}_W}P_{\widetilde{d}}(\xi, \xi')\widetilde{h}_W^{\pi_{\emptyset}}(\xi') - \widetilde{h}_W^{\pi_{\emptyset}}(\xi)\right] = 0.
    \end{align*}
    The left-hand side is exactly the policy evaluation equation $(\ref{pol_ev_2})$. Indeed, since $R$ is closed, the constant $C_R$ in expression (\ref{local_bias_av}) cancels. 
    Hence, for every recurrent class $R$, the function $h^{\pi_W}_{\text{av}}$ defined by (\ref{local_bias_av}) satisfies equation (\ref{pol_ev_2}). 

    Now, consider the functions $g^{\pi_W}_{\text{int}}$ and $h^{\pi_W}_{\text{int}}$ that arise from the intervention costs. Since $W$ is topple-free under $\pi_W$, grains can only leave $W$ through an intervention. Let $(\eta_t)_{t \geq 1}$ denote the process induced by policy $\pi_W$. Let
    \[
    I_T = \sum\limits_{t = 1}^{T} \mathbbm{1}\{d_W(\eta_t) \neq 0\}
    \]
    denote the number of interventions up to time $T$ and let
    \[
    J_T = \sum\limits_{t = 1}^{T} \sum\limits_{x \in W}[\eta_{t+1}(x) - (\eta_t^{d_W(\eta_t)})(x)]
    \]
    denote the number of grains that enter $W$ during the first $T$ transitions. Furthermore, let
    \[
    Z_t = \sum\limits_{x \in W} \eta_t(x), \qquad t \geq 1,
    \]
    denote the number of grains present in $W$ at time $t$. Since grains can only leave $W$ through an intervention and each intervention removes exactly one grain, we have
    \[
    Z_{T+1} - Z_1 = J_T - I_T, 
    \]
    or equivalently,
    \[
    I_T = J_T + Z_1 - Z_{T+1}.
    \]
    Since 
    \[
    0 \leq Z_t \leq \sum\limits_{x \in W} (\deg(x) - 1), \qquad \text{for all } t \geq 1,
    \]
    the quantity $Z_1 - Z_{T+1}$ is  bounded by the same constant for all $T \geq 1$. This implies
    \[
    \lim_{T \rightarrow \infty} \dfrac{1}{T}\mathbb{E}[I_T] = \lim_{T \rightarrow \infty} \dfrac{1}{T}\mathbb{E}[J_T].
    \]
    Thus, the long-run intervention rate equals the long-run rate at which grains enter $W$. An addition occurs directly in $W$ with probability $|W|/|V|$. If a new grain is added at a vertex $x \in V \setminus W$, then, by the coupled random walks representation, the probability it enters $W$ is given by $p_W(x)$.
    
    Thus, we obtain
    \begin{equation}\label{local_gain_int}
        g^{\pi_W}_{\text{int}}(\eta) = \dfrac{1}{|V|}\left(|W| + \sum\limits_{x \notin W}p_W(x)\right).
    \end{equation}
    Let $R$ denote a recurrent class of the Markov chain induced by $\pi_W$. Letting $\eta \in R$ and $\xi = \eta|_{V \setminus W}$, the expected future intervention cost generated by the sand grains that are present in $\eta$ is given by
    \begin{equation*}
        \left[\sum\limits_{x \in W}\eta(x) + \sum\limits_{x \notin W}p_W(x)\xi(x)\right]\beta.
    \end{equation*}
    Newly arriving sand grains contribute, in expectation, exactly $g^{\pi_W}_{\text{int}}$ per time step. This suggests that the intervention component of the bias is given by
    \begin{equation}\label{local_bias_int}
        h^{\pi_W}_{\text{int}}(\eta) = \sum\limits_{x \in W}\eta(x) + \sum\limits_{x \notin W}p_W(x)\eta(x) + C_R, \qquad \eta \in R.
    \end{equation}
    We verify this expression using the policy evaluation equation (\ref{pol_ev_2}). 

    Let 
    \begin{equation}\label{local_bias_int_help}
    H_{\text{int}}^{\pi_W}(\eta) = \sum\limits_{x \in W} \eta(x) + \sum\limits_{x \notin W}p_W(x) \eta(x), \qquad \eta \in \Omega.
    \end{equation}
    We define the function $q_W:V \cup \{s\} \rightarrow \mathbb{R}$ as
    \[
    q_W(x) = \begin{cases}
        1, &\text{if } x \in W,\\
        0, &\text{if } x = s, \\
        p_W(x), &\text{if } x \in V \setminus W.
    \end{cases}
    \]
    Then, expression (\ref{local_bias_int_help}) can be written as
    \[
    H_{\text{int}}^{\pi_W}(\eta) = \sum\limits_{x \in V}\eta(x)q_W(x), \qquad \eta \in \Omega.
    \]
    This implies
    \begin{equation}\label{local_bias_ex0}
        H^{\pi_W}_{\text{int}}(\eta^{d_W(\eta)}) - H^{\pi_W}_{\text{int}}(\eta) = -\mathbbm{1}\{d_W(\eta) \neq 0\}.
    \end{equation}
    For $\eta \in \Omega$ and $x \in V$, the definition of the stabilization process reads
    \[
    (\alpha_x \eta)^{\circ}(y) = \eta(y) + \mathbbm{1}_x(y) - \sum\limits_{u \in V}n_u(\eta, x) \Delta_{uy}, \qquad y \in V.
    \]
    Consequently,
    \begin{equation}\label{local_bias_ex1}
        H^{\pi_W}_{\text{int}}((\alpha_x \eta)^{\circ})- H^{\pi_W}_{\text{int}}(\eta) = q_W(x) - \sum\limits_{u \in V}n_u(\eta, x)(\Delta q_W)(u), \qquad \eta \in \Omega, \quad x \in V.
    \end{equation}
    Now, let $\eta \in R$. Since $W$ is topple-free under $\pi_W$, we have
    \[
    n_u(\eta^{d_W(\eta)},x) = 0, \qquad \text{for } u \in W, \quad x \in V.
    \]
    Moreover, since 
    \begin{align*}
        q_W(y) &= 1, \qquad y \in W, \\
        q_W(s) &= 0, \qquad \text{and }\\
        q_W(y) &= \dfrac{1}{\deg(y)}\sum\limits_{z \sim y}b_{yz}q_W(z), \qquad y \in V \setminus W,
    \end{align*}
    we have
    \[
    (\Delta q_W)(u) = \deg(u) q_W(u) - \sum\limits_{z \sim u}b_{uz}q_W(z) = 0, \qquad u \in V \setminus W.
    \]
    Therefore, using expression (\ref{local_bias_ex1}), we obtain 
    \begin{equation}\label{local_bias_ex2}
        H^{\pi_W}_{\text{int}}((\alpha_x \eta^{d_W(\eta)})^{\circ}) - H^{\pi_W}_{\text{int}}(\eta^{d_W(\eta)}) = q_W(x). 
    \end{equation}
    Combining expressions (\ref{local_bias_ex0}) and (\ref{local_bias_ex2}) now yields
    \[
    H^{\pi_W}_{\text{int}}((\alpha_x \eta^{d_W(\eta)})^{\circ}) - H^{\pi_W}_{\text{int}}(\eta) = q_W(x) - \mathbbm{1}\{d_W(\eta) \neq 0\}.
    \]
    Averaging over the addition vertex, we obtain
    \begin{equation}\label{local_bias_ex3}
    (P_{d_W}H^{\pi_W}_{\text{int}})(\eta) - H^{\pi_W}_{\text{int}}(\eta) = \dfrac{1}{|V|}\sum\limits_{x \in V}q_W(x) - \mathbbm{1}\{d_W(\eta) \neq 0\}.
    \end{equation}
    Since
    \[
    \dfrac{1}{|V|}\sum\limits_{x \in V}q_W(x) = \dfrac{1}{|V|}\left(|W| + \sum\limits_{x \notin W} p_W(x)\right) = g^{\pi_W}_{\text{int}},
    \]
    expression (\ref{local_bias_ex3}) reduces to 
    \[
    \mathbbm{1}\{d_W(\eta) \neq 0\} - g^{\pi_W}_{\text{int}} + (P_{d_W}H^{\pi_W}_{\text{int}})(\eta) - H^{\pi_W}_{\text{int}}(\eta) = 0,
    \]
    which is exactly the policy evaluation equation (\ref{pol_ev_2}). This shows that the intervention component of the bias is given by (\ref{local_bias_int}). 

    Substituting (\ref{local_gain_av}) and (\ref{local_gain_int}) into (\ref{local_gain_decomp}) and (\ref{local_bias_av}) and (\ref{local_bias_int}) into (\ref{local_bias_decomp}) yields (\ref{local_bias}) and (\ref{local_gain}), respectively. 
\end{proof}

Expression \eqref{local_gain} shows that the gain of the local intervention policy $\pi_W$ is the sum of the scaled gain of the never-intervene policy on the contracted graph and the long-run average intervention cost, determined by the probability that sand grains reach $W$ before hitting the sink. Similarly, up to an additive constant, the bias is the sum of the bias of the never-intervene policy on the contracted graph and an intervention term depending on the configuration, in which each grain is weighted by its probability of reaching $W$ before the sink. Using these, we now establish sufficient conditions for the optimality of a local intervention policy. 

\begin{thm}\label{general_optimal_policies}
    Let $\pi_W = (d_W)^{\infty}$ be any local intervention policy with intervention set $W \subseteq V$. Suppose that $W$ is topple-free under $\pi_W$. Let
    \[
        \underline{\beta}_W = \max_{\substack{a \in V \setminus W \\ \deg(a) \geq 2}} \dfrac{\phi_W(a)}{1-p_W(a)}, \qquad \overline{\beta}_W = \min_{u \in W: F_W(u)>0}\dfrac{1+\sum\limits_{v \in V \setminus W}b_{uv} \phi_W(v)}{F_W(u)},
    \]
    where
    \[
    F_W(u) = b_{us} + \sum\limits_{v \in V\setminus W}b_{uv}(1-p_W(v)).
    \]
    We adopt the conventions that the maximum over the empty set equals $0$, the minimum over the empty set equals $\infty$ and the ratio defining $\underline{\beta}_W$ is interpreted as $+\infty$ whenever $p_W(a) = 1$.
    If
    \[
    \underline{\beta}_W \leq \overline{\beta}_W,
    \]
    then $\pi_W$ is optimal for every $\beta \in [\underline{\beta}_W, \overline{\beta}_W]$. 
\end{thm}
\begin{proof}
    Define the function $q_W: V \cup \{0\} \rightarrow \mathbb{R}$ as
    \begin{align*}
        q_W(x) = \begin{cases}
            \beta, &\text{if } x \in W, \\
            \phi_W(x) + p_W(x)\beta, &\text{if } x \in V \setminus W, \\
            0, &\text{if } x = 0,
        \end{cases}
    \end{align*}
    where $0$ denotes the action corresponding to no intervention, whereas $s$ denotes the sink.
    Let \begin{equation}\label{bias_simple_gen}
        H_W(\eta) = \sum\limits_{y \in V}\eta(y)q_W(y).
    \end{equation}
    By expressions (\ref{local_bias}) and (\ref{h_pi0}), applied to the contracted graph $\widetilde{\mathcal{G}}_W$, the function $H_W$ coincides, up to an additive constant, with the bias $h^{\pi_W}$ on each recurrent class of the Markov chain induced by $\pi_W$. 

    Let $\eta \in \Omega$ and $x \in V$. By the definition of the stabilization process, we have
    \begin{equation}\label{stabilization}
    (\alpha_x \eta)^{\circ}(y) = \eta(y) + \mathbbm{1}_x(y) - \sum\limits_{u \in V}n_u(\eta, x) \Delta_{uy}, \qquad y \in V. 
    \end{equation}
    Using expressions (\ref{stabilization}) and (\ref{bias_simple_gen}), we obtain
    \begin{equation}\label{bias_diff_gen_1}
    H_W\left((\alpha_x \eta)^{\circ}\right) - H_W(\eta) = q_W(x) - \sum\limits_{u \in V}n_u(\eta, x)(\Delta q_W)(u).
    \end{equation}
    For any $u \in V \setminus W$, we have
    \[
    (\Delta \phi_W)(u) = 1, \qquad (\Delta p_W)(u) = 0.
    \]
    Hence, 
    \[
    (\Delta q_W)(u) = 1, \qquad u \in V \setminus W. 
    \]
    Substituting this into (\ref{bias_diff_gen_1}) yields
    \begin{equation}\label{ex_gen_0}
    H_W\left((\alpha_x \eta)^{\circ}\right) - H_W(\eta) = q_W(x) - \sum\limits_{u \in W}n_u(\eta, x)(\Delta q_W)(u) - \sum\limits_{u \in V \setminus W} n_u(\eta, x).
    \end{equation}
    Now, let $\eta \in R_{\pi_W}$. Since $W$ is topple-free under $\pi_W$, we have
    \[
    n_{u}(\eta^{d_W(\eta)}, x) = 0, \qquad u \in W, \quad x \in V.
    \]
    Therefore,
    \begin{equation}\label{ex_gen_1}
    H_W\left(\left(\alpha_x \eta^{d_W(\eta)}\right)^{\circ}\right) - H_W\left(\eta^{d_W(\eta)}\right) = q_W(x) - n\left(\eta^{d_W(\eta)}, x\right).
    \end{equation}
    Moreover, since $d_W(\eta) \in W$ whenever $d_W(\eta) \neq 0$ and $q_W(u) = \beta$ for every $u \in W$, expression (\ref{bias_simple_gen}) yields
    \begin{equation}\label{ex_gen_2}
    H_W\left(\eta^{d_W(\eta)}\right)- H_W(\eta) = -\beta \mathbbm{1}\{d_W(\eta) \neq 0\}.
    \end{equation}
    Combining (\ref{ex_gen_1}) and (\ref{ex_gen_2}) now yields
    \[
    H_W\left(\left(\alpha_x\eta^{d_W(\eta)}\right)^{\circ}\right) - H_W(\eta) = q_W(x) - n\left(\eta^{d_W(\eta)}, x \right) - \beta \mathbbm{1}\{d_W(\eta) \neq 0\}. 
    \]
    Consequently,
    \[
    n\left(\eta^{d_W(\eta)}, x\right) + \beta \mathbbm{1}\{d_W(\eta) \neq 0\} = q_W(x) - \left[H_W\left(\left(\alpha_x \eta^{d_W(\eta)}\right)^{\circ}\right) - H_W(\eta)\right].
    \]
    Using the definition of the cost function (\ref{cost_function}), we obtain, for every $\eta$ belonging to a recurrent class of the Markov chain induced by $\pi_W$,
    \begin{align*}
        c_{d_W}(\eta) &= \mathbb{E}\left[N\left(\eta^{d_W(\eta)}\right)\right] + \beta \mathbbm{1}\{d_W(\eta) \neq 0\} \\
        &= \dfrac{1}{|V|} \sum\limits_{x \in V} q_W(x) - \dfrac{1}{|V|} \sum\limits_{x \in V} \left[H_W\left(\left(\alpha_x \eta^{d_W(\eta)}\right)^{\circ}\right) - H_W(\eta)\right] \\
        &= \dfrac{1}{|V|}\sum\limits_{x \in V} q_W(x) + H_W(\eta) - \left(P_{d_W}H_W\right)(\eta).
    \end{align*}
    Moreover,
    \begin{align*}
        \dfrac{1}{|V|}\sum_{x \in V} q_W(x) &= \dfrac{1}{|V|}\sum\limits_{x \in V \setminus W} \phi_W(x) + \dfrac{1}{|V|}\left(|W| + \sum\limits_{x \in V \setminus W}p_W(x) \right)\beta \\
        &= \dfrac{|V|-|W|}{|V|}\widetilde{g}_W^{\pi_{\emptyset}} + \dfrac{1}{|V|} \left(|W| + \sum\limits_{x \in V \setminus W}p_W(x)\right)\beta \\
        &= g^{\pi_W},
    \end{align*}
    where the second equality follows from Theorem \ref{gain_bias_never_intervene}, applied to $\widetilde{\mathcal{G}}_W$, and the final equality follows from (\ref{local_gain}). Hence,
    \begin{equation}\label{policy_evaluation_W}
        c_{d_W}(\eta) = g^{\pi_W} + H_W(\eta) - (P_{d_W}H_W)(\eta)
    \end{equation}
    on every recurrent class of the Markov chain induced by $\pi_W$. 

    Now, let $\pi = d^{\infty}$ denote an arbitrary stationary deterministic policy. For $\eta \in \Omega$, $a \in A(\eta)$, and $x \in V$, expression (\ref{bias_simple_gen}) implies
    \[
    H_W(\eta^a) - H_W(\eta) = -q_W(a).
    \]
    Combining this with expression (\ref{ex_gen_0}), we obtain
    \begin{align*}
    H_W\left(\left(\alpha_x \eta^a\right)^{\circ}\right) - H_W(\eta) &= q_W(x) - \sum\limits_{u \in W}n_u(\eta^a, x)(\Delta q_W)(u) - \sum\limits_{u \in V \setminus W} n_u(\eta^a, x) - q_W(a) \\
    &= q_W(x) - n(\eta^a, x) + \sum\limits_{u \in W}n_u(\eta^a, x)[1-(\Delta q_W)(u)] - q_W(a).
    \end{align*}
    Rearranging gives 
    \begin{align}
    \nonumber n(\eta^a, x) + \beta \mathbbm{1}\{a \neq 0\} = &q_W(x) - [H_W((\alpha_x \eta^a)^{\circ}) - H_W(\eta)]\\
    \label{Bellman_decomp_W}&+ \beta \mathbbm{1}\{a \neq 0\} - q_W(a) + \sum\limits_{u \in W}n_u(\eta^a, x)[1-(\Delta q_W)(u)].
    \end{align}
    We now show that the final two terms in (\ref{Bellman_decomp_W}) are nonnegative when $\beta \in [\underline{\beta}_W, \overline{\beta}_W]$. 

    First, for $ a = d(\eta)$, we have
    \[
    \beta \mathbbm{1}\{a \neq 0\} - q_W(a) = \begin{cases}
        0, &\text{if } a = 0 \text{ or } a \in W, \\[1mm]
        \beta(1-p_W(a)) - \phi_W(a), &\text{if } a \in V \setminus W.
    \end{cases}
    \]
    This implies that, whenever
    \[
    \beta \geq \underline{\beta}_W = \max_{a \in V \setminus W} \dfrac{\phi_W(a)}{1-p_W(a)},
    \]
    we have
    \[
    \beta \mathbbm{1}\{a \neq 0\} - q_W(a) \geq 0.
    \]
    Next, let $u \in W$. Since $q_W(u) = \beta$, we have
    \begin{align*}
        (\Delta q_W)(u) &= \deg(u) \beta - \sum\limits_{v \in W}b_{uv}\beta - \sum\limits_{v \in V \setminus W} b_{uv}(\phi_W(v) + p_W(v)\beta) \\
        &= \beta \left[b_{us} + \sum\limits_{v \in V \setminus W} b_{uv}(1-p_W(v))\right] - \sum\limits_{v \in V \setminus W}b_{uv}\phi_W(v)
    \end{align*}
    Consequently,
    \begin{align*}
        1-(\Delta q_W)(u) &= 1+ \sum\limits_{v \in V \setminus W}b_{uv} \phi_W(v) - \beta\left[b_{us} + \sum\limits_{v \in V \setminus W}b_{uv}(1-p_W(v))\right]
    \end{align*}
    It follows that 
    \[
    1-(\Delta q_W)(u) \geq 0, \qquad u \in W,
    \]
    whenever
    \[
    \beta \leq \overline{\beta}_W = \min_{u \in W: F_W(u)>0}\dfrac{1+\sum\limits_{v \in V \setminus W}b_{uv} \phi_W(v)}{F_W(u)}.
    \]
   Thus, if $\beta \in [\underline{\beta}_W, \overline{\beta}_W]$, then expression (\ref{Bellman_decomp_W}) yields
   \[
   n(\eta^a, x) + \beta \mathbbm{1}\{a \neq 0\} \geq q_W(x) - \left[H_W\left((\alpha_x \eta^a)^{\circ}\right) - H_W(\eta)\right].
   \]
   Averaging over $x \in V$, we now obtain
   \begin{equation}\label{policy_evaluation_other}
   c_d(\eta) \geq \dfrac{1}{|V|}\sum\limits_{x \in V}q_W(x) + H_W(\eta) - (P_dH_W)(\eta) = g^{\pi_W} + H_W(\eta) - (P_dH_W)(\eta),
   \end{equation}
   for every $\eta \in \Omega$ and $a \in A(\eta)$. 

   To finalize the proof, let $(\eta_t)_{t \geq 1}$ denote the Markov chain induced by the arbitrary stationary deterministic policy $\pi = d^{\infty}$. Summing inequality (\ref{policy_evaluation_other}) over $t = 1, \ldots, T$ and taking expectations yields
   \[
   \mathbbm{E}_{\eta}^{\pi}\left[\sum\limits_{t = 1}^T c_d(\eta_t)\right] \geq Tg^{\pi_W} + H_W(\eta) - \mathbbm{E}_{\eta}^{\pi}[H_W(\eta_{T+1})].
   \]
   Dividing by $T$ and letting $T \rightarrow \infty$, using the fact that $H_W$ is bounded, now gives
   \[
   g^{\pi}(\eta) \geq g^{\pi_W} \qquad \eta \in \Omega.
   \]
   Thus, we conclude that the policy $\pi_W$ is optimal.   
\end{proof}

\section{Applications to specific graph classes}
\label{Applications_to_specific_graph_classes}

In this section, we apply the general results obtained in the preceding sections to four different graph classes: the one-dimensional lattice, the two-dimensional lattice, the complete graph and the star graph. Moreover, we discuss for each of these classes whether and under what conditions topple-freeness holds. For the one-dimensional lattice, the complete graph and the star graph, we characterize optimal policies for each value of the intervention cost.

\subsection{The one-dimensional lattice}
\label{The_one-dimensional_lattice}

In this subsection, we apply the general results on local intervention policies with topple-free intervention sets to characterize the optimal policies on the one-dimensional lattice for all values of $\beta$. 

Consider the one-dimensional lattice $\mathcal{G} = (V \cup \{s\}, E)$ with vertex set $V = \{1, \ldots, L\}$, where the boundary vertices $1$ and $L$ are connected to the sink by a single edge. The Green's function is given by \cite[Theorem 3]{CHUNG2000191}
\begin{equation}\label{Green}
    G(x,y) = \begin{cases}
        x(L+1-y)/(L+1), &\text{if } x \leq y, \\
        y(L+1-x)/(L+1), &\text{if } y < x.
    \end{cases}
\end{equation}
Using the definition (\ref{def_phi}), the function $\phi: V \rightarrow \mathbb{R}$ can be expressed as
\begin{equation}\label{phi_1D}
    \phi(x) = \sum_{y = 1}^{L}G(x,y) = \dfrac{x(L+1-x)}{2} = \dfrac{L+1}{2}G(x,x), \quad x \in V.
\end{equation}
Since $\phi(x)$ is a concave quadratic function of $x$, it is maximized at the central vertex when $L$ is odd and at the two central vertices when $L$ is even. Corollary \ref{gain_bias_opt_never_intervene_1D} provides explicit expressions for the gain, bias and optimality threshold of the never-intervene policy on the one-dimensional lattice. 

\begin{cor}\label{gain_bias_opt_never_intervene_1D}
    The gain of the never-intervene policy $\pi_{\emptyset} = (d_{\emptyset})^{\infty}$ is given by
    \begin{equation}\label{g_pi0_1D}
        g^{\pi_{\emptyset}}(\eta) = \dfrac{1}{12}(L+1)(L+2), \qquad \eta \in \Omega.
    \end{equation}
    Furthermore, there exists a constant $C \in \mathbb{R}$ such that 
    \begin{equation}\label{h_pi0_1D}
        h^{\pi_{\emptyset}}(\eta) = \sum\limits_{x = 1}^L (\eta(x) - 1) \phi(x) + C, \qquad \eta \in \Omega.
    \end{equation}
    The policy $\pi_{\emptyset}$ is optimal if and only if
    \begin{equation*}
        \beta \geq \begin{cases}
            (L+1)^2/8, &\text{if } L \text{ is odd}, \\
            L(L+2)/8, &\text{if } L \text{ is even}.
        \end{cases}
    \end{equation*}
\end{cor}
\begin{proof}
    The expressions for the gain and the bias follow immediately from Theorem \ref{gain_bias_never_intervene} and expression (\ref{phi_1D}). The optimality result follows immediately from Theorem \ref{nothing_optimal} and expression (\ref{phi_1D}). 
\end{proof}
We now consider local intervention policies with nested intervention sets that are central intervals of the lattice. For each $k \in \{0, \ldots, \lfloor L/2 \rfloor\}$, let
\[
W_k = \{k+1, \ldots, L-k\}
\]
and let $\pi_k = (d_k)^{\infty}$ be any local intervention policy with intervention set $W_k$. When $L$ is even, we have $W_{L/2} = \emptyset$, so $\pi_{L/2}$ coincides with the never-intervene policy. 

To describe the recurrent classes of the Markov chains induced by these policies, we introduce the configurations $\bar{\eta}_j \in \Omega$, for $j = 0, \ldots, \lfloor L/2 \rfloor$, defined as
\[
\bar{\eta}_j(i) = \begin{cases}
    1, &\text{if } i \in \{1, \ldots, j\} \cup \{L-j+1, \ldots, L\}, \\
    0, &\text{otherwise.}
\end{cases}
\]

The following theorem characterizes the unique recurrent class of the Markov chain induced by each policy $\pi_k$ for $k \in \{1, \ldots, \lfloor L/2 \rfloor\}$. The case $k = 0$ is treated separately in Proposition \ref{prop_always_intervene}.

\begin{thm}\label{rec_class_pi_k}
    For $k \in \{1, \ldots, \lfloor L/2 \rfloor\}$, let $\pi_k = (d_k)^{\infty}$ be any local intervention policy with intervention set $W_k$. 
    Then the Markov chain induced by $\pi_k$ has a unique recurrent class $R_k$, where
    \begin{equation}\label{def_Rk}
        R_k =
        \left\{
        \eta \in \{0,1\}^L :
        \sum_{i=1}^{k}\eta(i)\ge k-1,\;
        \sum_{i=L-k+1}^{L}\eta(i)\ge k-1,\;
        \sum_{i=k+1}^{L-k}\eta(i)\le 1
        \right\}
        \setminus\{\bar{\eta}_{k-1}\},
    \end{equation}
    for $k = 1, \ldots, \lfloor (L-1)/2 \rfloor$ and, when $L$ is even,
    \begin{equation}\label{def_Rl/2}
        R_{L/2} = \{\eta \in \{0,1\}^L: \sum\limits_{i=1}^L \eta(i) \geq L-1\}.
    \end{equation}
\end{thm}
\begin{proof}
    If $L$ is even and $k = L/2$, then $W_k = \emptyset$, so $\pi_k$ is the never-intervene policy. The uncontrolled sandpile dynamics on a one-dimensional lattice has a unique recurrent class consisting of stable configurations with at most one empty vertex \cite[Section 2]{Redig_survey}. This establishes expression (\ref{def_Rl/2}).

    Now, let $k \in \{1, \ldots, \lfloor (L-1)/2 \rfloor\}$.
    We first show that every recurrent configuration has at most one grain in the central strip $\{k+1, \ldots, L-k\}$. Whenever this strip contains at least one grain, the policy removes one of these grains before the addition of the new sand grain. Suppose that a configuration has at most one grain in the central strip. Then the strip is empty immediately after the intervention prescribed by the policy. The subsequent addition of a new grain followed by stabilization can create at most one grain in the central strip. It follows that the set
    \begin{equation}\label{rec_class_1}
    \{\eta \in \Omega: \sum\limits_{i = k+1}^{L-k} \eta(i) \leq 1\}
    \end{equation}
    is closed under the dynamics induced by $\pi_k$. 
    
    Now suppose that a configuration contains more than one grain in the central strip. At every transition under policy $\pi_k$, one such grain is removed. The subsequent addition and stabilization can create at most one new grain in the strip. Hence, the number of grains in the central strip cannot increase. Moreover, if the configuration outside the central strip is not completely filled, there is a positive probability that no new grain ends up in the central strip during the subsequent addition and stabilization, so that the number of central grains decreases by one. If the configuration outside the strip is completely filled, a new grain necessarily appears in the central strip during the next transition. However, there is a positive probability that a vertex outside the central strip becomes empty, after which there is again a positive probability that the number of grains in the central strip decreases at the following transition. Thus, from every configuration containing more than one grain in the central strip, there is a finite path of positive probability to a configuration with one fewer grain in the central strip. Iterating this argument shows that the set (\ref{rec_class_1}) can be reached with positive probability in finitely many steps. Since this set is closed, every recurrent class of the Markov chain induced by $\pi_k$ must be contained in (\ref{rec_class_1}).   

    Since the recurrent class is contained in (\ref{rec_class_1}), the central strip is empty immediately after an intervention. Moreover, the policy never removes grains from the two outside strips. Therefore, during the subsequent addition and stabilization, the dynamics restricted to each outside strip coincide with those of an uncontrolled sandpile model on a lattice of length $k$. Hence, the restriction of any recurrent configuration to either outside strip must belong to the recurrent class of the uncontrolled sandpile model on a lattice of length $k$, which is given by \cite{Redig_survey}
    \[
    \{\xi \in \{0,1\}^k: \sum\limits_{i=1}^k \xi(i) \geq k-1\}.
    \]
    It follows that every recurrent configuration has at most one empty vertex in each outside strip. Hence, every recurrent class is contained in
    \begin{equation}\label{rec_class_2}
    \left\{
        \eta \in \{0,1\}^L :
        \sum_{i=1}^{k}\eta(i)\ge k-1,\;
        \sum_{i=L-k+1}^{L}\eta(i)\ge k-1,\;
        \sum_{i=k+1}^{L-k}\eta(i)\le 1
        \right\}.
    \end{equation}
    We now show that from each configuration in $(\ref{rec_class_2})$, there exists a finite path of positive probability to the configuration $\bar{\eta}_k$. Let $\eta$ denote a configuration in the set $(\ref{rec_class_2})$. Suppose first that at least one of the outside strips contains an empty vertex. We then successively add grains at the empty vertices in the outside strips. Before each addition, the policy removes a grain from the central strip whenever one is present. Hence, after at most two transitions, both outside strips are completely filled and the central strip is empty, so that the system has reached $\bar{\eta}_k$.
    
    Now, suppose that both outside strips are already completely filled. If the central strip is empty, the configuration already equals $\bar{\eta}_k$. Otherwise, the central strip contains exactly one grain. The policy first removes this grain. Then, we add a grain at an arbitrary vertex in one of the outside strips. The resulting avalanche creates one empty vertex in that outside strip and one grain in the central strip. At the next transition, the policy removes the grain from the central strip, after which we add a grain at the newly created empty vertex. This yields $\bar{\eta}_k$. 

    Since each of these sequences of additions occurs with positive probability, it follows that $\bar{\eta}_k$ can be reached from every configuration in (\ref{rec_class_2}). Therefore, every recurrent class contains $\bar{\eta}_k$, which implies that the Markov chain induced by $\pi_k$ has a unique recurrent class.

    We proceed to show that the set of configurations accessible from $\bar{\eta}_k$ is precisely $R_k$, as defined in (\ref{def_Rk}). To this end, we show that every configuration in (\ref{rec_class_2}) can be reached from $\bar{\eta}_k$, except for the configuration $\bar{\eta}_{k-1}$. Let $\eta$ be a configuration in (\ref{def_Rk}). If $\eta = \bar{\eta}_k$, the statement is obvious. Suppose that $\eta \neq \bar{\eta}_k$. We distinguish seven cases.

    \paragraph{Case I:}
    Suppose that the outside strips in $\eta$ are completely filled. Since $\eta \neq \bar{\eta}_k$, there must be a grain in the central strip. In this case, $\eta$ can be reached from $\bar{\eta}_k$ by simply adding a grain at the corresponding vertex in the central strip. 

    \paragraph{Case II:}
    Suppose that $\eta$ has one empty vertex in either the left or the right outside strip, but not at a vertex adjacent to the central strip, and that there is no grain in the central strip. Without loss of generality, suppose that the empty vertex is $i \in \{1, \ldots, k-1\}$, and thus lies in the left outside strip. Starting from $\bar{\eta}_k$ we first add a grain at vertex $1$. After stabilization, this creates an empty vertex at $k$ and a grain at $k+1$. The policy then removes the grain at $k+1$, after which we add a grain at vertex $k-i$. This moves the empty vertex to $i$, yielding the configuration $\eta$.

    \paragraph{Case III:}
    Suppose that $\eta$ has one empty vertex at either the left or right neighbor of the central strip, and that there is no grain in the central strip. Without loss of generality, suppose that the empty vertex is $k$. Starting from $\bar{\eta}_k$, we again first add a grain at vertex $1$, creating an empty vertex at $k$ and a grain at $k+1$. After the policy removes the grain at $k+1$, we add a grain at an arbitrary vertex in the right outside strip. After stabilization, this creates an empty vertex in the right outside strip and a grain at $L-k$. The policy then removes the grain at $L-k$, after which we add a grain at the newly created empty vertex in the right outside strip. This results in the configuration $\eta$.

    \paragraph{Case IV:}
    Suppose that $\eta$ has one empty vertex in either the left or the right outside strip and one grain in the central strip. Without loss of generality, suppose that the empty vertex is $i \in \{1, \ldots, k\}$, and thus lies in the left outside strip. Starting from configuration $\bar{\eta}_k$, we first add a grain at vertex $k+1-i$. After stabilization, this creates an empty vertex at $i$ and a grain at vertex $k+1$. If the grain in the central strip of $\eta$ is located at $k+1$, we have now reached $\eta$. Otherwise, the policy first removes the grain at $k+1$, after which we add the grain at the vertex in the central strip occupied in $\eta$, yielding configuration $\eta$.

    \paragraph{Case V:}
    Suppose that $\eta$ has an empty vertex in both the left and the right outside strip, that at least one of these vertices is not adjacent to the central strip, and that there are no grains in the central strip. Without loss of generality, suppose that the empty vertices are $i \in \{1, \ldots, k\}$ and $j \in \{L-k+2, \ldots, L\}$, so that the empty vertex in the right outside strip is not adjacent to the central strip. Starting from configuration $\bar{\eta}_k$, we first add a grain at vertex $k+1-i$. After stabilization, this creates an empty vertex at $i$ and a grain at $k+1$. The policy then removes the grain at $k+1$. Next, we add a grain at vertex $L$, which creates an empty vertex at $L-k+1$ and a grain at $L-k$. After the policy removes the grain at $L-k$, we add a grain at vertex $2L+2-j-k$. After stabilization, the empty vertex has moved from $L-k+1$ to $j$, yielding the configuration $\eta$.

    \paragraph{Case VI:}
    Suppose that $\eta$ has an empty vertex in both the left and the right outside strip and that there is one grain in the central strip. Let the empty vertices be $i \in \{1, \ldots, k\}$ and $j \in \{L-k+1, \ldots, L\}$. Starting from configuration $\bar{\eta}_k$, we first add a grain at vertex $k+1-i$. After stabilization, this creates an empty vertex at $i$ and a grain at vertex $k+1$. The policy then removes the grain at $k+1$, after which we add a grain at vertex $2L-k-j+1$. After stabilization, this creates an empty vertex at $j$ and a grain at $L-k$. If $L-k$ is the occupied vertex in the central strip of $\eta$, we have reached $\eta$. Otherwise, the policy removes the grain at $L-k$, after which we add a grain at the vertex in the central strip that is occupied in $\eta$, yielding the configuration $\eta$. 

    \paragraph{Case VII:}
    Finally, consider the configuration $\bar{\eta}_{k-1}$. We show that this configuration cannot be reached in a single step from any configuration in (\ref{rec_class_2}) under the dynamics induced by $\pi_k$. Since the recurrent class $R_k$ is contained in (\ref{rec_class_2}), it then follows that $\bar{\eta}_{k-1}$ cannot be reached in a finite number of steps from $R_k$. Let $\eta$ denote a configuration in (\ref{rec_class_2}) from which a transition to $\bar{\eta}_{k-1}$ is possible, and let $i \in V$ denote the corresponding addition location. 
    
    Suppose first that $i \in \{k+1, \ldots, L-k\}$. Since the central strip is empty immediately after the intervention prescribed in $\eta$, addition at $i$ and subsequent stabilization leaves a grain at $i$. Therefore, the new configuration cannot be $\bar{\eta}_{k-1}$. 
    
    Now, suppose that $i \in \{1, \ldots, k\}$. For vertex $k$ to be empty after the addition at $i$ and subsequent stabilization, there are only two possibilities: either $i = 1$ and the left outside strip is completely filled, or $i \in \{2, \ldots, k\}$ and the left outside strip has one empty vertex at $i-1$. In either case, addition at $i$ and stabilization also creates a grain at vertex $k+1$. Hence, the resulting configuration cannot be $\bar{\eta}_{k-1}$. The same argument applies when $i \in \{L-k+1, \ldots, L\}$. Thus, the configuration $\bar{\eta}_{k-1}$ is not reachable from any configuration in (\ref{rec_class_2}).   

    We now conclude that $R_k$ is the unique recurrent class of policy $\pi_k$. We have shown that $\bar{\eta}_k$ belongs to every recurrent class, which implies that the recurrent class is unique. Furthermore, we have shown that this recurrent class is contained in (\ref{rec_class_2}). Cases I-VI show that each configuration in (\ref{rec_class_2}) except $\bar{\eta}_{k-1}$ is reachable from $\bar{\eta}_k$, whereas Case VII shows that $\bar{\eta}_{k-1}$ is not reachable. Thus, the unique recurrent class of policy $\pi_k$ is precisely $R_k$. 
\end{proof}

The above theorem establishes that for $k \in \{1, \ldots, \lfloor L/2 \rfloor\}$, any local intervention policy $\pi_k = (d_k)^{\infty}$ with intervention set $W_k$ induces a Markov chain with a unique recurrent class. The following example shows that for $k = 0$, this need not be the case. In Proposition \ref{prop_always_intervene}, however, we show that it is possible to construct an always-intervene policy with a unique recurrent class and that always-intervene policies with multiple recurrent classes cannot be optimal. 
\paragraph{Example}
Let $L \geq 4$ and consider the always-intervene policy $\pi^{\text{bnd}}_0 = (d^{\text{bnd}}_0)^{\infty}$, where 
\[
d^{\text{bnd}}_0(\eta) = \begin{cases}
    1, &\text{if } \eta(1) = 1, \\
    L, &\text{if } \eta(1) = 0 \text{ and } \eta(L) = 1, \\
    \min\{i \in \{1, \ldots, L\}|\eta(i) = 1\}, &\text{if } \eta(1) = \eta(L) = 0 \text{ and } \sum\limits_{i = 2}^{L-1}\eta(i) > 0, \\
    0, &\text{otherwise,}
\end{cases} \qquad \eta \in \Omega.
\]
We show that there exists a recurrent class that is contained in the set
\[
C = \{\eta \in \{0, 1\}^L: \sum\limits_{i=1}^L \eta(i) = 2, \sum\limits_{i=2}^{L-1} \eta(i) \geq 1\}.
\]
It suffices to show that $C$ is closed under the $\pi_0$-dynamics. Observe that any configuration in the set $C$ contains at least one grain at an interior vertex. Since each configuration contains exactly two grains and the policy favors removal from the boundary vertices, there is exactly one grain at an interior vertex $j$ after an intervention from any starting configuration in $C$. Suppose that the new grain is added at a vertex $k \neq j$. In that case, the new configuration is contained in $C$. If $k = j$, then the vertex $j$ topples and the new configuration contains grains at $j-1$ and $j+1$. Since $L \geq 4$, at least one of these vertices is an interior vertex. Hence, the resulting configuration belongs to $C$. This shows that $C$ is closed and therefore contains at least one recurrent class of the Markov chain induced by $\pi_0$. Since $C$ is disjoint from the recurrent class $R_0$, it follows that the Markov chain induced by $\pi_0$ has at least two recurrent classes.
\begin{prop}
    \label{prop_always_intervene}
    For the case $k = 0$, the following statements hold:
    \begin{enumerate}[leftmargin=*]
        \item Under every always-intervene policy $\pi_0 = (d_0)^{\infty}$, the set
        \begin{equation}\label{def_R0}
        R_0 = \{\eta \in \{0,1\}^L: \sum\limits_{i = 1}^L \eta(i) = 1\}
        \end{equation}
        is a recurrent class of the induced Markov chain.
        \item There exists an always-intervene policy $\pi_0 = (d_0)^{\infty}$ under which $R_0$ is the unique recurrent class of the induced Markov chain.
        \item Let $\pi_0$ and $\pi_0'$ be always-intervene policies whose induced Markov chains have a unique recurrent class and multiple recurrent classes, respectively. Then,
        \[
        g^{\pi_0}(\eta) = \beta \leq g^{\pi_0'}(\eta), \qquad \eta \in \Omega.
        \]
        The inequality above is strict for any $\eta$ that belongs to a recurrent class $R \neq R_0$ of $\pi_0'$. Thus, $\pi_0'$ is not optimal for any $\beta \geq 0$. 
    \end{enumerate}
\end{prop}
\begin{proof}\quad\\
    \begin{enumerate}[leftmargin=*]
        \item It suffices to show that the set $R_0$ is closed and irreducible under the $\pi_0$-dynamics. For every $\eta \in R_0$, the policy removes the unique grain present in $\eta$. The subsequent addition results in another configuration in $R_0$. This implies that $R_0$ is closed. Moreover, every configuration in $R_0$ is reached during a single transition with probability $1/L$, regardless of the starting configuration in $R_0$. Hence, the set $R_0$ is irreducible. Thus, $R_0$ is a recurrent class of the Markov chain induced by any always-intervene policy. 
        \item Consider the always-intervene policy $\pi^{\text{left}}_0 = (d^{\text{left}}_0)^{\infty}$, where
        \[
        d^{\text{left}}_0(\eta) = \begin{cases}
            0, &\text{if } \eta(i) = 0 \text{ for all } i \in \{1, \ldots, L\}, \\
            \min\{i \in \{1, \ldots, L\}: \eta(i) = 1\}, &\text{otherwise.}
        \end{cases}
        \]
        We show that under this policy, $R_0$ can be reached in a finite number of transitions from any starting configuration. It then follows from the first statement that $R_0$ is the unique recurrent class. 

        Let $\eta$ denote a configuration containing at least two grains and let $j$ denote the position of its rightmost grain. We specify a sequence of transitions that occurs with positive probability and reaches the set $R_0$. Since the policy $\pi^{\text{left}}_0$ always removes the leftmost grain, the grain at $j$ remains after the intervention. Suppose that $j < L$. Then, we add the next grain at vertex $j+1$. By assumption, this vertex is empty, so the resulting configuration is stable and the rightmost occupied vertex is now $j+1$. Repeating this procedure, we reach a configuration whose rightmost grain is at vertex $L$. During these transitions, the total number of grains remains constant. 
        Once the rightmost grain is at vertex $L$, we add a grain at $L$ after the intervention. After the resulting avalanche,  exactly one grain is lost to the sink. Hence, the total number of grains decreases by one. This procedure can be repeated until a configuration is reached with exactly one grain, which happens in a finite number of transitions. Since each transition has positive probability of occurring, the set $R_0$ is accessible from every configuration. Thus, invoking the first statement, $R_0$ is the unique recurrent class of the Markov chain induced by $\pi_0$.
        \item Within the recurrent class $R_0$, one sand grain is removed during each transition and no topplings occur. Therefore,
        \[
        g^{\pi_0}(\eta) = g^{\pi_0'}(\eta) = \beta \qquad \text{for all } \eta \in R_0.
        \]
        Let $R \neq R_0$ denote a recurrent class of the Markov chain induced by $\pi_0'$. Every configuration in $R$ contains at least two grains. Therefore, during each transition, one grain is removed and the probability of a toppling is at least $1/L$. Therefore,
        \[
        g^{\pi_0'}(\eta) \geq \beta + \dfrac{1}{L} > \beta \qquad \eta \in R.
        \]
        By expression (\ref{av_cost_def_2}),  the gain at any configuration $\eta \in \Omega$ is a convex combination of the stationary average costs on each of the recurrent classes. Since each recurrent class of the Markov chain induced by $\pi_0'$ has stationary average cost at least $\beta$ and $R_0$ is the unique recurrent class of the Markov chain induced by $\pi_0$,  we obtain
        \[
        g^{\pi_0'}(\eta) \geq \beta = g^{\pi_0}(\eta), \qquad \eta \in \Omega,
        \]
        where the inequality is strict on every recurrent class $R \neq R_0$. Thus, the policy $\pi_0'$ is not optimal for any $\beta \geq 0$.
    \end{enumerate}
\end{proof}
Given a configuration $\eta \in \Omega$ and $k = 1, \ldots, \lfloor L/2 \rfloor$, we define the configurations
\[
\eta^{(k)}_{\ell}(i) = \eta(k+1-i), \qquad \eta^{(k)}_r(i) = \eta(L-k+i), \qquad i = 1, \ldots, k.
\]
These are the restrictions of $\eta$ to the leftmost and rightmost $k$ vertices, respectively, indexed from the interior towards the boundary. For $k = 0$, we regard $\eta_{\ell}^{(0)}$ and $\eta_r^{(0)}$ as configurations on the empty lattice. 
The following theorem gives the gain and the bias of the policies $\pi_k$, $k = 0, \ldots, (L-1)/2$. 
\begin{thm}\label{recurrent_class_pi_k}
    For each $k \in \{1, \ldots, \lfloor (L-1)/2 \rfloor\}$, let $\pi_k = (d_k)^{\infty}$ be any local intervention policy with intervention set $W_k$. Let $\pi_0 = (d_0)^{\infty}$ denote an always-intervene policy whose induced Markov chain has the unique recurrent class $R_0$, defined in \eqref{def_R0}. Then, for every $k \in \{0, \ldots, \lfloor (L-1)/2 \rfloor \}$, the gain is given by
    \begin{equation}\label{gain_pi_k}
        g^{\pi_k}(\eta) = \dfrac{k(k+1)(k+2)}{6L} + \dfrac{L-k}{L}\beta, \qquad \eta \in \Omega.
    \end{equation}
    Moreover, for every $k \in \{0, \ldots, \lfloor (L-1)/2 \rfloor \}$, there exists a constant $C_k \in \mathbb{R}$ such that for every $\eta \in R_k$, the bias is given by
    \begin{equation}\label{bias_pi_k}
        h^{\pi_k}(\eta) = \beta \sum\limits_{x = k+1}^{L-k} \eta(x) + \sum\limits_{i = 1}^k \left[\dfrac{i(k+1-i)}{2} + \beta \left(1-\dfrac{i}{k+1}\right)\right](\eta^{(k)}_{\ell}(i) + \eta^{(k)}_r(i)) + C_k.
    \end{equation}
    For odd $L$ and $k = (L-1)/2$, expression (\ref{bias_pi_k}) holds for every $\eta \in \Omega$.
\end{thm}
\begin{proof}
    We first show that for each $k \in \{0, \ldots, \lfloor (L-1)/2 \rfloor \}$, the set $W_k$ is topple-free under the corresponding policy $\pi_k$. The result then follows by applying Theorem \ref{gain_bias_pi_W_topple_free}. For $k = 0$, the Markov chain induced by $\pi_0$ has a unique recurrent class given by $R_0$. Since each configuration in $R_0$ contains exactly one grain, which is removed by the policy before the next addition, no topplings occur during transitions from $R_0$. Therefore, the set $W_0 = V$ is topple-free under $\pi_0$.  Now, consider $k \in \{1, \ldots, \lfloor (L-1)/2 \rfloor \}$. By the definition (\ref{def_Rk}) of the recurrent class $R_k$ of the Markov chain induced by $\pi_k$, each $\eta \in R_k$ has at most one sand grain in the set $W_k$. Thus, immediately after the intervention by the policy, the set $W_k$ is empty. Therefore, it is impossible for any vertex $x \in W_k$ to topple during the stabilization following the next addition. Hence, the set $W_k$ is topple-free under the policy $\pi_k$ for each $k \in \{0, \ldots, \lfloor (L-1)/2 \rfloor \}$.

    We now invoke Theorem \ref{gain_bias_pi_W_topple_free}. Let $\widetilde{g}^{\pi_{\emptyset}}_{W_k}$ and $\widetilde{h}^{\pi_{\emptyset}}_{W_k}$ denote the gain and bias of the never-intervene policy on the contracted graph $\widetilde{\mathcal{G}}_{W_k}$. For $k = 0$, the contracted graph $\widetilde{\mathcal{G}}_{W_k}$ has no non-sink vertices. We adopt the conventions $\widetilde{g}_V^{\pi_{\emptyset}} = 0$ and $\widetilde{h}_V^{\pi_{\emptyset}} = 0$. For any $k \in \{0, \ldots, \lfloor (L-1)/2 \rfloor \}$, we have \cite[Section 12.2]{grinstead1997introduction}
    \begin{equation}\label{hitting_prob_1D}
    p_{W_k}(k+1-i) = p_{W_k}(L-k+i) = 1 - \dfrac{i}{k+1}, \qquad i = 1, \ldots, k.
    \end{equation}
    Substituting this and $|V| = L$ and $|W| = L-2k$ into expression (\ref{local_gain}), we obtain
    \[
    g^{\pi_k}(\eta) = \dfrac{2k}{L}\widetilde{g}^{\pi_{\emptyset}}_{W_k} + \dfrac{L-k}{L}\beta.
    \]
    For $k = 0$, this gives
    \[
    g^{\pi_0}_L(\eta) = \beta, \qquad \eta \in \Omega,
    \]
    in agreement with expression (\ref{gain_pi_k}).
    For $k \geq 1$, this graph consists of two one-dimensional lattices of length $k$, each connected to the contracted sink. Since a new grain is added to either component with probability $1/2$, expression (\ref{g_pi0_1D}) yields (\ref{gain_pi_k}). 

    Furthermore, the general bias formula (\ref{h_pi0}) shows that the bias on the contracted graph is the sum of the biases on its two components, up to an additive constant. Substituting expressions (\ref{phi_1D}), (\ref{h_pi0_1D}) and (\ref{hitting_prob_1D})  into (\ref{local_bias}) yields expression (\ref{bias_pi_k}). 

    Finally, suppose that $L$ is odd and $k = (L-1)/2$, so that $W_k = \{z\}$. For every $\eta \in \Omega$, the vertex $z$ is empty immediately after applying the intervention prescribed by $\pi_k$. If the next grain is added at $z$, no toppling occurs. Otherwise, only vertices in the connected component of occupied vertices containing the addition site in the post-decision configuration can topple, so $z$ cannot topple. This implies that $z$ never topples under policy $\pi_k$, even during transitions from transient configurations. Therefore, expression (\ref{bias_pi_k}) satisfies the policy evaluation equation on the entire configuration space. Since the chain has a unique recurrent class, expression (\ref{bias_pi_k}) holds for every $\eta \in \Omega$.
\end{proof}
Finally, Theorem \ref{1D_optimal_policies} identifies an optimal local intervention policy for each value of $\beta$. 

\begin{thm}\label{1D_optimal_policies}
    Let $K = \left\lfloor \dfrac{L-1}{2}\right\rfloor$.
    Define 
    \[
        \beta_k = \dfrac{k(k+1)}{2}, \qquad k = 0, \ldots, K, 
    \]
    and
    \[
        \beta_{K+1} = 
        \begin{cases}
            (L+1)^2/8, &\text{if } L \text{ is odd},\\
            L(L+2)/8, &\text{if } L \text{ is even}.
        \end{cases}
    \]
    For each $k = 1, \ldots, K$, let $\pi_k = (d_k)^{\infty}$ be any local intervention policy with intervention set $W_k$. Let $\pi_0$ be an always-intervene policy whose induced Markov chain has unique recurrent class $R_0$. Then, for each $k = 0, \ldots, K$, 
    \[
    \pi_k \text{ is optimal} \quad \Longleftrightarrow \quad \beta \in [\beta_k, \beta_{k+1}].
    \]
    Moreover,
    \[
    \pi_{\emptyset} \text{ is optimal} \quad \Longleftrightarrow \quad \beta \geq \beta_{K+1}.
    \]
\end{thm}
\begin{proof}
    We first prove sufficiency. In the proof of Theorem \ref{recurrent_class_pi_k}, we established that $W_k$ is topple-free under $\pi_k$ for every $k = 0, \ldots, K$. We can therefore invoke Theorem \ref{general_optimal_policies}. 
    For $k \geq 1$, the contracted graph $\widetilde{\mathcal{G}}_{W_k}$ consists of two paths of length $k$ that are both connected to the contracted sink. For $i = 1, \ldots, k$, we have, by \eqref{phi_1D},
    \[
    \phi_{W_k}(k+1-i) = \phi_{W_k}(L-k+i) = \dfrac{i(k+1-i)}{2},
    \]
    Combining this expression with \eqref{hitting_prob_1D}, we obtain
    \[
    \underline{\beta}_{W_k} = \max_{i \in \{1, \ldots, k\}}\dfrac{(k+1-i)(k+1)}{2} = \dfrac{k(k+1)}{2} = \beta_k.
    \]
    For $k = 0$, the maximum is taken over the empty set, so $\underline{\beta}_{W_0} = \beta_0 = 0$. 

    We proceed to evaluate the upper bound. First, suppose that $|W_k| \geq 2$. If $u \in W_k$ is not an endpoint of $W_k$, we have
    \[
    F_{W_k}(u) = 0.
    \]
    If $u$ is an endpoint of $W_k$, then,  using \eqref{hitting_prob_1D} with $i = 1$, we obtain
    \[
    F_{W_k}(u) = \dfrac{1}{k+1}
    \]
    and
    \[
    1 + \sum_{v \in V \setminus W_k} b_{uv} \phi_{W_k}(v) = \dfrac{k+2}{2}.
    \]
    Therefore, 
    \[
    \overline{\beta}_{W_k} = \dfrac{(k+1)(k+2)}{2} = \beta_{k+1}.
    \]
    For even $L$, this includes the case $k = K = L/2-1$. 
    
    Suppose now that $W_k$ is a singleton. This implies that $L$ is odd and $k = K = (L-1)/2$. Let the central vertex be denoted by $z$. We have
    \[
    F_{W_k}(z) = \dfrac{2}{K+1} \qquad \text{and } 1 + \sum_{v \in V \setminus W_K} b_{zv} \phi_{W_K}(v) = K + 1. 
    \]
    Thus,
    \[
    \overline{\beta}_{W_k} = \dfrac{(L+1)^2}{8} = \beta_{K+1}.
    \]
    Theorem \ref{general_optimal_policies} now implies that $\pi_k$ is optimal whenever $\beta \in [\beta_k, \beta_{k+1}]$. 

    For the never-intervene policy, we take $W = \emptyset$ and again invoke Theorem \ref{general_optimal_policies}. Since $p_{\emptyset} = 0$ and $\phi_{\emptyset} = \phi$, we obtain
    \[
    \underline{\beta}_{\emptyset} = \max_{x \in V}\phi(x) = \max_{x \in \{1, \ldots, L\}}\dfrac{x(L+1-x)}{2} = \beta_{K+1}
    \]
    and 
    \[
    \overline{\beta}_{\emptyset} = \infty.
    \]
    Thus, $\pi_{\emptyset}$ is optimal whenever $\beta \geq\beta_{K+1}$. 

    We next prove necessity. Suppose that $\beta < \beta_k$, $k = 1, \ldots, K$. 
    By expression (\ref{gain_pi_k}), we have
    \[
        g^{\pi_k}(\eta) - g^{\pi_{k-1}}(\eta) = \dfrac{\beta_k - \beta}{L} > 0, \qquad k = 1, \ldots, K.
    \]
    Hence, if $\beta < \beta_k$, the policy $\pi_{k-1}$ has strictly smaller gain than $\pi_k$ for all $k = 1, \ldots, K$. 

    Now, suppose that $\beta > \beta_{k+1}$, $k = 0, \ldots, K-1$. 
    By expression (\ref{gain_pi_k}), we have
    \[
    g^{\pi_k}(\eta) - g^{\pi_{k+1}}(\eta) = \dfrac{\beta - \beta_{k+1}}{L}.
    \]
    Hence, if $\beta > \beta_{k+1}$, the policy $\pi_{k+1}$ has strictly smaller gain than $\pi_k$. 

    Finally, using \eqref{g_pi0_1D}, we obtain both for odd and even $L$,
    \[
    g^{\pi_K}(\eta) - g^{\pi_{\emptyset}}(\eta) = \dfrac{(L-K)(\beta - \beta_{K+1})}{L}.
    \]
    It follows that $\pi_K$ is not optimal when $\beta > \beta_{K+1}$ and $\pi_{\emptyset}$ is not optimal when $\beta < \beta_{K+1}$. This establishes necessity for each policy $\pi_k$, $k = 0, \ldots, K$. 
\end{proof}

\subsection{The two-dimensional lattice}
\label{The_two_dimensional_lattice}

\paragraph{The two-dimensional lattice}
Consider the two-dimensional square lattice $\mathcal{G} = (V \cup \{s\}, E)$ with vertex set $V = \{1, \ldots, L\}^2$, where every edge leaving $V$ is connected to the sink $s$. The Green's function $G: V \times V \rightarrow \mathbb{R}$ is given by the following spectral representation \citep[Section 5]{CHUNG2000191} (note that the formula in \citep{CHUNG2000191} differs by a 
factor of $4$ from ours, since we consider the continuous-time random walk that takes 
jumps at rate $4$, rather than the discrete-time walk)
\begin{equation}\label{Green_2D}
    G(x, y) = \dfrac{4}{(L+1)^2}\sum\limits_{k=1}^L\sum\limits_{\ell = 1}^L \dfrac{\sin\left(\dfrac{\pi k x_1}{L+1}\right)\sin\left(\dfrac{\pi k y_1}{L+1}\right) \sin\left(\dfrac{\pi \ell x_2}{L+1}\right)\sin \left(\dfrac{\pi \ell y_2}{L+1}\right)}{4 - 2\cos \left(\dfrac{\pi k }{L+1}\right) - 2\cos \left(\dfrac{\pi \ell}{L+1}\right)}.
\end{equation}
Combining (\ref{Green_2D}) with definition (\ref{def_phi}) yields an explicit expression for the function $\phi$, from which the gain and bias functions follow via (\ref{g_pi0}) and (\ref{h_pi0}). To determine the critical threshold of $\beta$ for which the never-intervene policy is optimal, it remains to identify the maximizers of $\phi$. 
\begin{cor}\label{cor_threshold_never_intervene_2D}
    Consider the two-dimensional square lattice $\mathcal{G} = (V \cup \{s\}, E)$, where $V = \{1, \ldots, L\}^2$ and every edge leaving $V$ is connected to the sink $s$. The policy $\pi_{\emptyset} = (d_{\emptyset})^{\infty}$, where $d_{\emptyset}(\eta) = 0$ for all $\eta \in \Omega$, is optimal if and only if
    \[
        \beta \geq
            \begin{cases}
            \phi\left(\left(\frac{L+1}{2},\frac{L+1}{2}\right)\right),
            &\text{if }L\text{ is odd},\\[2mm]
            \phi\left(\left(\frac{L}{2},\frac{L}{2}\right)\right),
            &\text{if }L\text{ is even},
            \end{cases}
    \]
    where, in the even case,
    \[
        \phi\left(\left(\frac{L}{2},\frac{L}{2}\right)\right)
        =
        \phi\left(\left(\frac{L}{2},\frac{L}{2}+1\right)\right)
        =
        \phi\left(\left(\frac{L}{2}+1,\frac{L}{2}\right)\right)
        =
        \phi\left(\left(\frac{L}{2}+1,\frac{L}{2}+1\right)\right).
    \]
\end{cor}
\begin{proof}
    To prove the statement, we make use of the probabilistic interpretation of the function $\phi$ as the expected first hitting time to the sink of a continuous-time random walk $(X_t)_{t\geq 0}$ with jump rate 1 along each edge,  as given by \eqref{prob_interpret}. We write $X_t = (X_t^{(1)}, X_t^{(2)})$. Observe that the two coordinates are independent one-dimensional continuous-time simple random walks. Let
    \[
    \tau^{(k)} = \inf\{t \geq 0: X^{(k)}_t \notin \{1, \ldots, L\}\}, \qquad k = 1, 2.
    \]
    Then
    \[
    \tau_s = \min\{\tau^{(1)}, \tau^{(2)}\}.
    \]
    Therefore,
    \[
    \phi(x_1, x_2) = \mathbb{E}_{(x_1, x_2)}[\tau_s] = \mathbb{E}_{(x_1, x_2)}[\min\{\tau^{(1)}, \tau^{(2)}\}] = \int\limits_0^{\infty}\mathbb{P}_{x_1}(\tau^{(1)} > t)\mathbb{P}_{x_2}(\tau^{(2)} > t)dt.
    \]
    Since $\mathbb{P}_i(\tau^{(k)} > t)$, $k = 1, 2$, $i = 1, \ldots, L$, is nondecreasing as $i$ moves toward the center, it follows that $\phi$ is maximized at the central vertex $((L+1)/2, (L+1)/2)$ if $L$ is odd or at the four central vertices $(L/2, L/2)$, $(L/2, L/2 + 1)$, $(L/2 + 1, L/2)$ and $(L/2 + 1, L/2 + 1)$ if $L$ is even. The result now follows from Theorem \ref{nothing_optimal}.
\end{proof}

By a numerical investigation, we found that intervention sets on the two-dimensional lattice need not be topple-free under the corresponding local intervention policies. Therefore, the results under the topple-freeness condition do not apply in general in this case.

\subsection{The complete graph}
\label{The_complete_graph}

Consider the complete graph $\mathcal{G} = (V \cup \{s\}, E)$, where $V = \{1, \ldots, n\}$. Suppose that each vertex $x \in V$ is connected to the sink. In this case, the restriction of the graph Laplacian of $\mathcal{G}$ to $V \times V$ is given by
\begin{equation*}
    \Delta_{xy} = \begin{cases}
        n, &\text{if } x = y, \\
        -1, &\text{if } x \neq y,
    \end{cases}
\end{equation*}
that is,
\begin{equation*}
    \Delta = (n+1)I_n - J_n,
\end{equation*}
where $I_n$ denotes the $n \times n$ identity matrix and $J_n$ the $n \times n$ all-ones matrix. This implies that the Green's function is given by
\begin{equation*}
    G = \Delta^{-1} = \dfrac{1}{n+1}(I_n + J_n),
\end{equation*}
or equivalently,
\begin{equation*}
    G(x,y) = \begin{cases}
        \dfrac{2}{n+1}, &\text{if } x = y, \\[1.5ex]
        \dfrac{1}{n+1}, &\text{if } x \neq y.
    \end{cases}
\end{equation*}
Using definition (\ref{def_phi}), we now obtain
\begin{equation}\label{phi_complete_graph}
    \phi(x) = \sum\limits_{y \in V} G(x,y) = \dfrac{2}{n+1} + (n-1)\dfrac{1}{n+1} = 1, \quad \text{for all } x \in V.
\end{equation}
The gain and bias of the policy $\pi_{\emptyset}$ and the optimality threshold for never intervening in case of the complete graph are stated in Corollary \ref{gain_bias_opt_never_intervene_complete}.
\begin{cor}\label{gain_bias_opt_never_intervene_complete}
    The gain of the never-intervene policy $\pi_{\emptyset} = (d_{\emptyset})^{\infty}$ is given by
    \[
        g^{\pi_{\emptyset}}(\eta) = 1, \qquad \eta \in \Omega.
    \]
    Furthermore, there exists a constant $C \in \mathbb{R}$ such that
    \begin{align*}
        h^{\pi_{\emptyset}}(\eta) = \sum\limits_{x \in V}(\eta(x) - n +1) + C, \qquad \eta \in \Omega.
    \end{align*}
    The policy $\pi_{\emptyset}$ is optimal if and only if $\beta \geq 1$.
\end{cor}
\begin{proof}
    The expressions for the gain and the bias follow immediately from Theorem \ref{gain_bias_never_intervene} and expression (\ref{phi_complete_graph}). The optimality result follows immediately from Theorem \ref{nothing_optimal} and expression (\ref{phi_complete_graph}). 
\end{proof}

\begin{cor}\label{complete_graph_always_intervene}
    Let $\pi_V$ denote an always-intervene policy. Its gain is given by
    \[
    g^{\pi_V}(\eta) = \beta, \qquad \eta \in \Omega
    \]
    and it is optimal if and only if $\beta \leq 1$.
\end{cor}
\begin{proof}
    Under $\pi_V$ the mass cannot increase from a nonempty configuration. Hence, within each recurrent class of the Markov chain induced by $\pi_V$, the mass must remain constant. Since every toppling dissipates one sand grain on the complete graph, no topplings can occur within a recurrent class. Therefore, $V$ is topple-free under $\pi_V$. Moreover, one intervention occurs at every time step. This yields.
    \[
    g^{\pi_V}(\eta) = \beta, \qquad \eta \in \Omega.
    \]
    We now apply Theorem \ref{general_optimal_policies} with $W = V$. For the lower bound, we obtain
    \[
    \underline{\beta}_V = 0.
    \]
    For the upper bound, we have
    \[
    F_V(u) = b_{us} = 1, \qquad u \in V
    \]
    and hence,
    \[
    \overline{\beta}_V = 1.
    \]
    Thus, $\pi_V$ is optimal whenever $\beta \leq 1$. 

    Conversely, it follows from Corollary \ref{gain_bias_opt_never_intervene_complete} that the never-intervene policy has gain $g^{\pi_{\emptyset}}(\eta) = 1$ for every $\eta \in \Omega$. Therefore, the always-intervene policy $\pi_V$ is not optimal if $\beta > 1$. 
\end{proof}

Corollaries \ref{complete_graph_always_intervene} and \ref{gain_bias_opt_never_intervene_complete} provide a complete characterization of the optimality regimes in terms of $\beta$. 

\begin{cor}
For $\beta \geq 1$, the never-intervene policy is optimal, whereas for $\beta \leq 1$, every always-intervene policy is optimal. In particular, when $\beta = 1$, both types of policies are optimal.
\end{cor}

We next investigate topple-freeness on the complete graph. The following Lemma shows that singleton intervention sets need not be topple-free in the complete graph. We checked numerically by an exhaustive search, that in the complete graph on $\{ 1, \dots, n \} \cup \{ s \}$ each vertex is topple-free when $1 \le n \le 6$. We show below that this is no longer true when $n \ge 7$.

\begin{lem}
    If $n \geq 7$, no vertex $x \in V$ is topple-free under the single-vertex intervention policy $\pi_x$. 
\end{lem}
\begin{proof}
    Let $x \in V$ and let $\eta \in \mathcal{R}_{\pi_x}$. By sequentially adding particles to noncritical vertices in $V \setminus \{x\}$, we can reach a configuration $\overline{\eta}$ that satisfies
    \[
    \overline{\eta}(y) = \begin{cases}
        k, &\text{if } y = x, \\
        n-1, &\text{otherwise}
    \end{cases}
    \]
    for some $k = 0, 1, \ldots, n-1$. Since $\eta$ is a recurrent configuration and this sequence of additions occurs with positive probability, it follows that $\overline{\eta}$ is a recurrent configuration as well. We now show that, from every such configuration $\overline{\eta}$, there exists a finite sequence of transitions, occurring with positive probability, during which $x$ topples. 
    
    First, suppose that $k \geq 2$. Choose any vertex $z \in V \setminus \{x\}$. Before adding a new grain at $z$, the policy reduces the height at $x$ to $k-1$. Since all vertices in $V \setminus \{x\}$ have critical height, adding a particle at $z$ causes every vertex in $V \setminus \{x\}$ to topple. As a result, the vertex $x$ receives $n-1$ particles and reaches height $k + n - 2 \geq n$. This implies that $x$ topples. 

    Now, consider $k \in \{0,1\}$. In either case, after application of $d_x$, the height at $x$ equals $0$. Choose any vertex $z \in V \setminus \{x\}$ and add a sand grain at $z$. This causes every vertex in $V \setminus \{x\}$ to topple. Let the resulting configuration be denoted by $\eta'$. We have
    \[
    \eta'(y) = \begin{cases}
        n-1, &\text{if } y = x, \\
        n-2, &\text{if } y = z, \\
        n-3, &\text{otherwise.}
    \end{cases}
    \]
    Now, let $y_1, y_2$ denote distinct vertices in $V \setminus \{x, z\}$. Starting from $\eta'$, perform four successive transitions with addition sites $z, y_1, y_1, y_2$ respectively. Note that none of these additions causes a toppling and before each addition one particle is removed from $x$. Since $n \geq 7$, the resulting configuration $\eta''$ is given by
    \[
    \eta''(u) = \begin{cases}
        n-1, &\text{if } u \in \{z, y_1\}, \\
        n-2, &\text{if } u = y_2, \\
        n-5, &\text{if } u = x, \\
        n-3, &\text{otherwise.}
    \end{cases}
    \]
    Finally, add another sand grain at $z$. Before this addition, the policy reduces the height at $x$ to $n-6$. Observe that this addition causes every vertex in $V \setminus \{x\}$ to topple. It follows that $x$ reaches height $2n-7$. Since $n \geq 7$ implies $2n-7 \geq n$, the vertex $x$ topples. Since the constructed sequence of transitions starts in a recurrent configuration and occurs with positive probability, $x$ is not topple-free under $\pi_x$. 
\end{proof}

\subsection{The star graph}
\label{The_star_graph}

Consider the star graph $\mathcal{G} = (V \cup \{s\}, E)$, where $V = \{z, 1, \ldots, n\}$ and $z$ denotes the center. Suppose that only the center is connected to the sink. Ordering the vertices as $z, 1, \ldots, n$, the restriction of the graph Laplacian of $\mathcal{G}$ to $V \times V$ is given by
\[
\Delta_{xy} = \begin{cases}
    n+1, &\text{if } x = y = z, \\
    -1, &\text{if } x = z \neq y, \quad \text{or } y = z \neq x, \\
    1, &\text{if } x = y \neq z, \\
    0, &\text{otherwise.}
\end{cases}
\]
Hence, the Green's function is given by
\[
G(x,y) = \Delta^{-1}_{xy} = \begin{cases}
    2, &\text{if } x = y \neq z, \\
    1, &\text{otherwise.}
\end{cases}
\]
Invoking definition (\ref{def_phi}), we now find
\begin{equation}\label{phi_star_graph}
\phi(x) = \sum\limits_{y \in V} G(x,y) = \begin{cases}
    n+1, &\text{if } x = z, \\
    n+2, &\text{otherwise.}
\end{cases}
\end{equation}
Using this result, we now provide the gain and bias of the policy $\pi_{\emptyset}$, together with the threshold for never intervening, for the star graph.

\begin{cor}\label{gain_bias_opt_never_intervene_star}
    The gain of the never-intervene policy $\pi_{\emptyset} = (d_{\emptyset})^{\infty}$ is given by
    \begin{equation}\label{gain_never_intervene_star}
        g^{\pi_{\emptyset}}(\eta) = n+2 - \dfrac{1}{n+1}, \qquad \eta \in \Omega.
    \end{equation}
    Furthermore, there exists a constant $C \in \mathbb{R}$ such that
    \begin{equation}\label{bias_never_intervene_star}
        h^{\pi_{\emptyset}}(\eta) = (\eta(z) - n)(n+1) + C, \qquad \eta \in \Omega.
    \end{equation}
    The policy $\pi_{\emptyset}$ is optimal if and only if $\beta \geq n+1$.
\end{cor}
\begin{proof}
    The expressions for the gain and the bias follow immediately from Theorem (\ref{gain_bias_never_intervene}) and expression (\ref{phi_star_graph}). The optimality result follows immediately from Theorem \ref{nothing_optimal} and expression (\ref{phi_star_graph}). 
\end{proof}
We next show that the center $z$ is topple-free under the single-vertex intervention policy $\pi_z$.
\begin{lem}\label{star_topple-free}
    The center $z$ is topple-free under the policy $\pi_z$. 
\end{lem}
\begin{proof}
    After applying the decision rule $d_z$, the height at $z$ is at most $n-1$. Since every leaf has degree 1, adding a grain at a leaf causes it to topple once and send one grain to $z$. Adding a grain at $z$ also 
    increases its height by one. Hence, during the addition and subsequent stabilization, the height at $z$ increases by exactly one. This implies that $z$ never topples under the policy $\pi_z$. 
\end{proof}
Since every leaf has height zero in every stable configuration, the policy $\pi_z$ coincides with the always-intervene policy $\pi_V$. The following lemma provides the gain of this policy.
\begin{lem}\label{value_always_intervene_star}
    The gain of the always-intervene policy $\pi_z = (d_z)^{\infty}$ is given by
    \begin{equation}\label{gain_always_intervene_star}
    g^{\pi_z}(\eta) = \dfrac{n}{n+1} + \beta, \qquad \eta \in \Omega.
    \end{equation}
\end{lem}
\begin{proof}

    In Lemma \ref{star_topple-free}, we established that $z$ is topple-free under the policy $\pi_z$. Hence, we can apply Theorem \ref{gain_bias_pi_W_topple_free} with $W = \{z\}$. On the contracted graph $\widetilde{\mathcal{G}}_z$, each of the $n$ leaves is connected to the sink by one edge. Hence, every addition causes exactly one toppling. This yields
    \[
    \widetilde{g}_z^{\pi_{\emptyset}} = 1.
    \]
    Furthermore,
    \[
    p_z(x) = 1, \qquad x \in V \setminus \{z\}.
    \]
    Applying Theorem \ref{gain_bias_pi_W_topple_free} now yields
    \[
    g^{\pi_z}(\eta) = \dfrac{n}{n+1} + \beta, \qquad \eta \in \Omega.
    \]
\end{proof}

We now establish a necessary and sufficient condition under which the always intervene policy $\pi_z$ is optimal.
\begin{cor}\label{always_intervene_star}
    The always-intervene policy $\pi_z$ is optimal if and only if 
    \[
    \beta \leq n+1.
    \]
\end{cor}
\begin{proof}
    By Lemma \ref{star_topple-free}, the intervention set $W = \{z\}$ is topple-free under $\pi_z$. Hence, we can apply Theorem \ref{general_optimal_policies}.
    Since each of the leaves has degree one, the maximum that defines the lower bound is taken over the empty set. This implies
    \[
    \underline{\beta}_z = 0.
    \]
    Moreover, we have for each leaf $x$,
    \[
    p_z(x) = 1, \qquad \text{and } \phi_z(x) = 1.
    \]
    Therefore, we obtain
    \[
    \overline{\beta}_z = n+1.
    \]
    Thus, the policy $\pi_z$ is optimal whenever $\beta \leq n+1$.

    To prove necessity, we compute, using expressions (\ref{gain_never_intervene_star}) and (\ref{gain_always_intervene_star}),
    \[
    g^{\pi_z}(\eta) - g^{\pi_{\emptyset}}(\eta) = \beta - (n+1).
    \]
    Thus, when $\beta > n+1$, the gain of the never-intervene policy is strictly smaller than that of $\pi_z$, which implies that $\pi_z$ is not optimal.
\end{proof}

Corollaries \ref{gain_bias_opt_never_intervene_star} and \ref{always_intervene_star} provide the following complete characterization of the optimal policies on the star graph. 

\begin{cor}
    The always-intervene policy is optimal for $\beta \leq n+1$, whereas the never-intervene policy is optimal for $\beta \geq n+1$. At $\beta = n+1$, both policies are optimal.
\end{cor}

\section{Discussion}
\label{Discussion}

In this paper, we formulated and analyzed the problem of controlling avalanche sizes through targeted removal of sand grains as an average-cost MDP that balances avalanche size against intervention cost. Within this formulation, we studied local intervention policies. We derived general structural results for such policies and sharpened these results under the condition of topple-freeness. In addition, for policies that satisfy this condition, we presented explicit expressions for the gain and bias and explicit bounds on $\beta$ that guarantee their optimality. We applied our general results to several specific graph classes, leading, among other things, to a complete characterization of the optimal policies on the one-dimensional lattice, the complete graph and the star graph for all values of $\beta$.

This work opens several interesting avenues for future research. First, the characterization of optimal policies on the one-dimensional lattice naturally raises the question whether similarly structured optimal policies exist for other graph classes. Our results show that the critical value of $\beta$ above which the never-intervene policy is optimal is determined by the maximum of the function $\phi$. On the one-dimensional lattice, as the intervention cost decreases, the optimal intervention sets expand to include increasingly many vertices with the largest values of $\phi$. It is therefore natural to investigate whether similarly structured local intervention policies are optimal on other graph classes, such as higher-dimensional lattices. Moreover, it would be interesting to study the dependence of the size and shape of the corresponding intervention sets on the intervention cost.


Second, since several results rely on the property of topple-freeness, it would be interesting to determine whether this condition can be characterized in terms of the topology of the underlying graph or the spanning trees associated with recurrent sandpile configurations. Such graph-theoretic conditions could facilitate the identification of topple-free intervention sets to which our results apply.

Third, a possible direction for future research is to study the control problem under a discounted-cost criterion. Such a criterion assigns greater relative weight to the costs incurred in the near future, where the relative importance of future costs is determined by the discount factor. It therefore allows us to investigate the trade-off between reducing immediate costs and accounting for the subsequent evolution of the system. This would interpolate between minimizing immediate avalanche costs, as studied by \citet{deJongh1}, and minimizing the long-run average cost, as considered in the present work. Depending on the value of the discount factor, the resulting optimal policies may differ substantially from those obtained under the average-cost criterion. The discounted-cost criterion would be particularly relevant when interventions are intended to mitigate imminent avalanches.


Fourth, when an analytical solution is unavailable, approximate dynamic programming or reinforcement learning methods could be used to construct approximately optimal policies. The exact structural results obtained in this work could guide the design of reduced state and policy spaces. One possible approach is to aggregate configurations into a smaller collection of representative states or to represent them through a set of carefully selected features.  

Finally, future research could consider alternative intervention mechanisms, such as controlling the addition site or temporarily modifying the toppling threshold. Different cost functions could also be studied. In this work, the cost function captures the trade-off between the immediate intervention cost and the expected avalanche size. Alternative objectives could instead penalize the probability that an avalanche exceeds a prescribed size or assign a nonlinear cost to large-scale avalanches. Such objectives may lead to policies that more strongly prioritize the prevention of catastrophic events.

%
%
%
%

\bigskip
\par\noindent
\textbf{Acknowledgments}
\par\noindent
MCdJ thanks AAJ and the Department of Mathematical Sciences at the University of Bath for their kind hospitality while conducting this research. MCdJ is also grateful to the Catharine van Tussenbroek Fund for financial support for her research visit to Bath.

\printbibliography

\end{document}